\RequirePackage{fix-cm}
\documentclass[smallextended]{svjour3}
\smartqed
\usepackage[utf8]{inputenc}
\usepackage[T1]{fontenc}
\usepackage{microtype}
\usepackage{mathptmx}
\usepackage{mathtools}
\usepackage[UKenglish]{babel}
\usepackage{csquotes}
\usepackage{graphicx}
\usepackage{amsmath}
\usepackage{amssymb}
\usepackage{multirow}
\usepackage{array}                                                  
\usepackage{booktabs}                                               
\usepackage{dcolumn}                                                
\usepackage{multirow}
\usepackage{tikz}
\usepackage{siunitx}
\usepackage[english,algoruled,linesnumbered,longend]{algorithm2e}
\usepackage[acronym]{glossaries}
\usepackage[sort,compress]{cite}
\DeclareMathOperator*{\argmin}{arg\,min}

\DeclareMathOperator{\vecop}{vec}
\DeclareMathOperator{\traceOp}{tr}
\DeclareMathOperator{\shrink}{shrink}

\DeclareMathOperator{\Rpart}{Re}
\DeclareMathOperator{\Ipart}{Im}
\newcommand{\Fro}{\mathrm{F}}
\newcommand{\Real}{\mathbb{R}}
\newcommand{\Complex}{\mathbb{C}}
\newcommand{\trace}[1]{\traceOp\left(#1\right)}
\newcommand{\set}[1]{\left\lbrace#1\right\rbrace}
\newcommand{\norm}[1]{\left\|#1\right\|}
\newcommand{\abs}[1]{\left|#1\right|}
\newcommand{\scprod}[2]{\langle#1,#2\rangle}
\newcommand{\coloneqq}{\mathrel{\mathop:}=}
\SetAlFnt{\normalsize}
\SetAlCapSkip{0.5ex}
\usetikzlibrary{%
  decorations.pathreplacing,%
  shapes,%
  plotmarks,%
  arrows,%
  topaths,%
  patterns,%
  external,%
  circuits.ee.IEC%
}
\begin{document}
% - Meta Data ------------------------------------------------------------------
\title{%
  Sparse $\ell_{1}$ Regularisation of Matrix Valued Models\\
  for Acoustic Source Characterisation
  \thanks{This work has partially been funded by the German Research Foundation
    (DFG) within the grant number SA~1502/5-1. This funding is thankfully
    acknowledged.}
}
%
% \subtitle{Revision:~\today}
\titlerunning{Sparse Acoustic Source Characterisation}
\author{%
  Laurent~Hoeltgen\and
  Michael~Breuß\and
  Gert~Herold\and
  Ennes~Sarradj
}
\authorrunning{Hoeltgen et al.}
\institute{%
  L.\ Hoeltgen\and M.\ Breuß\at
  Chair for Applied Mathematics,
  Brandenburg University of Technology Cottbus-Senftenberg\\
  Platz der Deutschen Einheit 1,
  03046 Cottbus, Germany\\
  \email{$\lbrace$hoeltgen, breuss$\rbrace$@b-tu.de}
  \and
  G. Herold\and E. Sarradj\at
  Chair for Technical Acoustics,
  Brandenburg University of Technology Cottbus-Senftenberg\\
  Siemens-Halske-Ring 14,
  03046 Cottbus, Germany\\
  \email{$\lbrace$herold, sarradj$\rbrace$@b-tu.de}
}
\date{Received: date / Accepted: date}
\newacronym{FBS}{FBS}{Forward Backward Splitting}
\newacronym{SVD}{SVD}{Singular Value Decomposition}
\newacronym{PDHG}{PDHG}{Primal Dual Hybrid Gradient}
\newacronym{SB}{SB}{Split Bregman}
\newacronym{BP}{BP}{Basis Pursuit}
\newacronym{GPGPU}{GPGPU}{General Purpose Computation on Graphics Processing
  Unit}
\newacronym{RMS}{RMS}{Root Mean Square}
\newacronym{CMF}{CMF}{Covariance Matrix Fitting}
\newacronym{CLEANSC}{Clean-SC}{Clean based on Source Coherence}
\newacronym{FFT}{FFT}{Fast Fourier Transform}
% ------------------------------------------------------------------------------
\maketitle
% ------------------------------------------------------------------------------
\begin{abstract}
      We present a strategy for the recovery of a sparse solution of a common
      problem in acoustic engineering, which is the reconstruction of sound
      source levels and locations applying microphone array measurements. The
      considered task bears similarities to the basis pursuit formalism but also
      relies on additional model assumptions that are challenging from a
      mathematical point of view. Our approach reformulates the original task as
      a convex optimisation model. The sought solution shall be a matrix with a
      certain desired structure. We enforce this structure through additional
      constraints. By combining popular splitting algorithms and matrix
      differential theory in a novel framework we obtain a numerically efficient
      strategy. Besides a thorough theoretical consideration we also provide an
      experimental setup that certifies the usability of our strategy. Finally,
      we also address practical issues, such as the handling of inaccuracies in
      the measurement and corruption of the given data. We provide a post
      processing step that is capable of yielding an almost perfect solution in
      such circumstances.
\end{abstract}
% 
% ------------------------------------------------------------------------------
\section{Introduction}
\label{sec:Introduction}
% ------------------------------------------------------------------------------
% 
In 2005, Osher \emph{et al.}~\cite{Osher2005} proposed an algorithm for the
iterative regularisation of inverse problems that was based on findings of
Bregman~\cite{Bregman1967}. They used this algorithm, nowadays called Bregman
iteration, for image restoration purposes such as denoising and deblurring.
Especially in combination with the Rudin-Osher-Fatemi (ROF) model for denoising
\cite{ROF92} they were able to produce excellent results. Their findings caused
a subsequent surge of interest in the Bregman iteration. Among the numerous
application fields, it has for example been used to solve the basis pursuit
problem \cite{Cai2009,Osher2008,Yin2007} and was later applied to medical
imaging problems \cite{Lin2006}. Further applications include deconvolution and
sparse reconstructions \cite{ZBBO09}, wavelet based denoising \cite{Xu2006} and
nonlinear inverse scale space methods \cite{Burger2006,Burger2005}. An important
adaptation of the Bregman iteration is the \gls{SB} method \cite{GO2009} and the
linearised Bregman approach \cite{Cai2009}. The \gls{SB} algorithm can be used
to solve $\ell_{1}$-regularised inverse problems in an efficient way. Its
benefits stem from the fact that differentiability is not a necessary
requirement on the underlying model.\par
\begin{figure}
      \centering
      \begin{tikzpicture}[>=latex]
            \draw[thick] (2,-0.3) node [circle, draw] {};
            \draw[thick] (3.5,0.8) node [circle, draw] {};
            \draw[thick] (3.6,1.8) node[circle, draw] {};
            \draw[thick] (2.8,-1.2) node[circle, draw] {};
            \draw[thick] (2.2, 2) node (M)
            [circle, draw,label=above right:$c_j$] {};
            \draw[thick] (2.5, 0.8) node (C) [circle, draw] {};
            \draw[thick] (-2, 1.8) node (C) [circle, draw, fill] {};
            \draw[thick] (-2.5, 1) node (S)
            [circle, fill, label=below left:$x_s$] {};
            \draw[thick] (-3., -1.5) node (T) [cross out,draw] {};
            \draw[thick] (-2.8, 2) node (T) [cross out, draw] {};
            \draw[thick] (-3.7, 1) node (T) [cross out, draw] {};
            \draw[thick] (-3.5, -1) node (T)
            [cross out, draw,label=below left:$x_i$] {};
            \draw[thick,->] (S) -- (M) node [midway, above] {$a_{sj}$};
            \draw[thick, dashed, ->] (T) -- (M) node 
            [near start, above] {$a_{ij}$};
            \draw (S.north west) node
            [above, text width=1cm] {sources};
            \draw (T.north west) node
            [above left,text width=2cm,text ragged left]
            {possible source locations };
            \draw (3,-0.5) node
            [above, text width=2cm, text centered] {microphone array};
      \end{tikzpicture}
      \caption{A generic set-up for acoustic source characterisation using a
        microphone array. Each source $s$ causes a sound pressure amplitude
        $x_{s}$ at the reference point. The sound pressure $c_{j}$ arriving at
        microphone $j$ is given by $c_{j}=a_{sj}x_{s}$. If multiple sources are
        present, then it follows from a superposition principle that we have
        $c_{j}=\sum_{i}a_{ij}x_{i}$.}
      \label{fig:setup}
\end{figure}
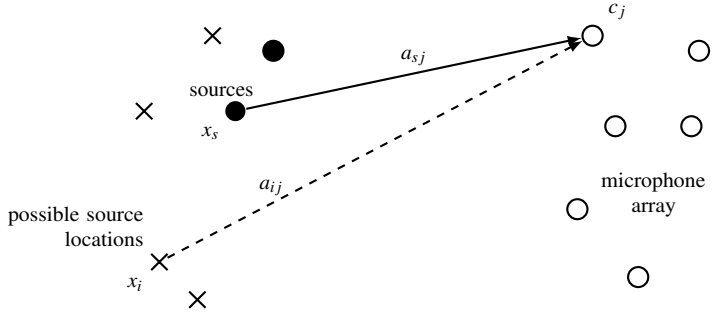

Our contribution is concerned with an application to acoustic source
characterisation using a microphone array. A microphone array comprises $n$
microphones at known locations. These microphones register the sound that is
emitted by a number of sources with unknown locations (see
Fig.~\ref{fig:setup}). The characterisation of these acoustic sources requires
the estimation of their location and strength.\par
We briefly describe the physical background behind the models that we discuss in
this work. The propagation of sound from a source position $x$ to a receiver at
position $y$ can be modelled by a Green's function. In the following we assume
that the source is always a monopole. In that case the sound pressure amplitude
at the receiver position for a given discrete frequency $\omega$ is defined by
\begin{equation}
      \label{eq:1}
      p \left( r, \omega \right) =  q_{0} \frac{1}{r}
      \exp\left( -\imath \omega \frac{r}{c_{0}} \right)
\end{equation}
with $\imath$ being the complex unit, $q_{0}$ being the source strength, and $r$
denoting the distance between the source and the receiver. Finally, the constant
$c_{0}$ denotes the speed of sound. The signals from any given source are
evaluated at a reference point $y_{0}$ at distance $r_{0}$ from the source:
\begin{equation}
      \label{eq:2}
      p_{0} \left( r_{0}, \omega \right) =  q_{0} \frac{1}{r_{0}}
      \exp\left( -\imath \omega \frac{r_{0}}{c_{0}} \right)
\end{equation}
Introducing the reference point into our model formulation helps us to eliminate
the source strength and leads to the following description
\begin{equation}
      \label{eq:3}
      p\left( r, \omega \right) = a p_{0} \left( r_{0}, \omega \right)
\end{equation}
with
\begin{equation}
      \label{eq:4}
      a \coloneqq \frac{r_{0}}{r}
      \exp\left( \imath \omega \frac{r_{0}-r}{c_{0}} \right)
\end{equation}
Equation~\eqref{eq:3} yields the sound pressure amplitude at a receiver
position depending on the sound pressure induced at a reference location by the
source.\par
The estimation approach that we follow here additionally assumes that the actual
locations of the sources are restricted to $m$ possible source locations. Since
a superposition principle holds in our model, we can account for multiple
sources by adding all contributions. Thus, the sound pressure at a microphone
$j$ is given by $\sum_i a_{ij} x_i$ where the sum is taken over all possible
source locations $i$. The coefficient $a_{ij}$ is defined in accordance with
\eqref{eq:4}:
\begin{equation}
      \label{eq:5}
      a_{ij} \coloneqq \frac{r_{0,j}}{r_{ij}}
      \exp\left(\imath \omega \frac{r_{0,j}-r_{ij}}{c_{0}} \right)
      \quad i = 1, \ldots, m,\ \text{and}\ j = 1, \ldots, n
\end{equation}
with $r_{ij}$ denoting the distance between sender $i$ and receiver $j$.\par
In a typical setting the number $m$ of possible source locations is much larger
than the number $n$ of microphones, but the number of actual sources is less
than $n$. Therefore, most of the $x_{i}$ in the sum are zero. The described
setup is also visualised in Fig.\ref{fig:setup}.\par
Gathering all possible coefficients $a_{ij}$ in a matrix $A$, and using a
vector $c \in \Complex^{n}$ to hold all microphone sound pressures yields
\begin{equation}
      \label{eq:propmodel}
      c=Ax,
\end{equation}
where $x \in \Complex^m$ is the sparse vector of source strengths. Using the
Hermitian cross-spectral matrix $C \coloneqq E[cc^\top]$ of microphone sound
pressures, we can reformulate \eqref{eq:propmodel} as
\begin{equation}
      \label{eq:6}
      C=A\underbrace{E\left[xx^\top\right]}_{\eqqcolon X} A^\top
\end{equation}  
where the operator $E$ denotes the expected value and where $X \in \Complex$
is the cross spectral matrix of source levels. This matrix is sparse,
Hermitian, and in case of uncorrelated source signals also diagonal. By
estimating $X$, the task of characterising the sources is solved. Let us
emphasise that this is an ill-posed inverse problem. We refer the reader to
\cite{HS2015a, HS2014, HSG2013} for further details.
% 
% ------------------------------------------------------------------------------
\paragraph{Our Contribution.}
\label{sec:Contribution}
% ------------------------------------------------------------------------------
% 
In this paper we discuss different models to characterise sound sources. Our
starting point is \eqref{eq:6}. Instead of solving this linear system of
equations directly we suggest to minimise the corresponding least squares
formulation and to augment our model with a sparsity inducing regulariser. Two
different approaches can be taken into account to respect the assumption that
the solution $X$ from \eqref{eq:6} should be a sparse diagonal matrix. We
present a formulation that adds harsh penalisation weights to off-diagonal
entries, and a second formulation that seeks optimal solutions over the space of
diagonal matrices. Both approaches use $\ell_{1}$ norms of the solution $X$ to
enhance sparsity.\par
Further, we propose numerical strategies for both models. Our basic algorithm is
an adapted version of the popular \gls{SB} method, which has been used
successfully for sparsity constrained problems in the past. It is well suited to
handle the non-differentiability of our cost function. For our second model we
additionally introduce advanced matrix valued differentiation concepts to
preserve the diagonal structure of the solution. We see the main novelties of
our work in these two ideas, namely the application of the \gls{SB} scheme on
the task at hand as well as the usage of matrix differential calculus to enforce
structural constraints in the solution.\par
Finally, we present a post processing strategy to handle setups with corrupted
data. The benefits of our approach are a robust behaviour with very sparse
solutions. Our solvers are tailored such that the number of non-zero entries in
their solutions coincides exactly with the number of sought sources.
% 
% ------------------------------------------------------------------------------
\paragraph{Structure of the Paper.}
\label{sec:Structure}
% ------------------------------------------------------------------------------
% 
In Section~\ref{sec:ProblemFormulation} the considered problem is discussed in
more detail. We elaborate on the occurring difficulties and propose a certain
number of models to overcome the mentioned difficulties in
Section~\ref{sec:analys-sugg-models}. Section~\ref{sec:SolvingStrategy} presents
our considered numerical approaches. The necessary mathematical preliminaries
are briefly outlined in Section~\ref{sec:Preliminaries} whereas
Sections~\ref{sec:BregmanFramework}-\ref{sec:OptSplitBregStruct} provide a
detailed presentation of our numerical schemes. Potential post processing
methods to improve the results further are established in
Section~\ref{sec:post-proc-steps}. Afterwards we provide a thorough evaluation
on real world data sets in Section~\ref{sec:NumericalEvaluation}. Finally, we
give a short conclusion with an outlook on challenges and future work in
Section~\ref{sec:ConclusionOutlook}.
% 
% ------------------------------------------------------------------------------
\section{Problem Formulation}
\label{sec:ProblemFormulation}
% ------------------------------------------------------------------------------
% 
The mathematical formulation of the considered problem is as follows. Given
matrices $A\in\Complex^{n,m}$, and $C\in\Complex^{n,n}$, where $C$ is also known
to be Hermitian, we seek a diagonal matrix $X\in\Complex^{m,m}$ with a sparse
set of non-zero entries along its main diagonal such that the equation
\begin{equation}
  \label{eq:7}
  AXA^{\top} = C
\end{equation}
holds. We remark that the previous equation can easily be rewritten as
\begin{equation}
  \label{eq:8}
  \left(A\otimes{}A\right)\vecop{(X)} = \vecop{(C)}
\end{equation}
where $\otimes$ represents the Kronecker matrix product and where the $\vecop$
operator stacks the entries from its argument column wise on top of each other.
We refer to \cite{PP2008,MN2007,M2000a} for a thorough presentation and
properties of these two operators. It follows that one could actually reduce the
initial problem in \eqref{eq:7} to solving the linear system given in
\eqref{eq:8}. However, the system matrix is dense and prohibitively large and
thus unsuited for most computational routines. If each entry in $A\otimes{}A$ is
stored in double precision, then we require $16n^{2}m^{2}$ byte of memory. For
the experimental setups that we consider in this paper this equals to about
185~GB of data. In addition, we cannot assert that all our model assumptions
will always be fulfilled if we simply solve the linear system. Indeed, in
practice the matrices $A$ and $C$ will stem from experimental setups and be
corrupted with noise. This means, even if we were to find a solution to
\eqref{eq:7}, we cannot expect it to be a real valued diagonal matrix with a
sparse set of non-zero entries along its main diagonal. Indeed, we have the
following result, which can also be found verbatim in \cite{MN2007} (Chapter~2,
Theorem~3)
\begin{proposition}
      \label{thm:1}
      A necessary and sufficient condition for the matrix equation $AXB = C$ to
      have a solution is that $AA^{\dagger} C B^{\dagger} B = C$, in which case
      the general solution is
      \begin{equation}
            \label{eq:9}
            X = A^{\dagger} C B^{\dagger} + Q - A^{\dagger} A Q B B^{\dagger}
      \end{equation}
      where $Q$ is an arbitrary matrix of appropriate order.
\end{proposition}
Here, $A^{\dagger}$ denotes the Moore-Penrose inverse of the matrix $A$.
Applying the previous proposition onto \eqref{eq:7} with $B=A^{\top}$ yields
conditions that assert the existence of a solution. Also, it becomes obvious
that the solution of the system may not necessarily have the sought
structure.\par
In order to overcome the aforementioned restrictions we suggest to consider the
following variational models
\begin{gather}
  \label{eq:10}
  \argmin_{X\in\Complex^{m,m}}
  \set{\frac{1}{2}\norm{AXA^{\top}-C}_{\Fro}^{2} + \lambda \norm{X}_{1}} \\
  \label{eq:11}
  \argmin_{X\in\Complex^{m,m}}
  \set{\frac{1}{2}\norm{AXA^{\top}-C}_{\Fro}^{2} + \norm{W\circ{}X}_{1}} \\
  \label{eq:12}
  \left\lbrace
  \begin{gathered}
    \argmin_{X\in\Complex^{m,m}}
    \set{\frac{1}{2}\norm{AXA^{\top}-C}_{\Fro}^{2} + \lambda \norm{X}_{1}} \\
    \text{under the constraint that $X$ is a diagonal matrix}
  \end{gathered}\right.
\end{gather}
where $\circ$ represents the Hadamard matrix product \cite{M2007d}. It simply
multiplies all matrix entries componentwise with each other. Also,
$\norm{\cdot}_{\Fro}$ represents the Frobenius norm for matrices
\cite{HJ1990,HJ1994}. In order to treat complex valued arguments, we
additionally use the following conventions:
\begin{align}
  \label{eq:13}
  \norm{X}_{1} &\coloneqq
  \sum_{i,j} \left( \abs{\Rpart (X_{ij})} + \abs{\Ipart (X_{ij})} \right) \\
  \label{eq:14}
  \norm{X}_{\Fro{}}^{2} &\coloneqq
  \sum_{i,j} \left( \Rpart (X_{ij})^{2} + \Ipart (X_{ij})^{2} \right)
\end{align}
Here, $\Rpart(z)$ and $\Ipart(z)$ represent the real and imaginary part of the
complex number $z$. Our choice is motivated by the fact that it will allow us a
fast parallel optimisation in the construction of our numerical schemes. An
in-depth analysis on the models is given in the forthcoming
Section~\ref{sec:analys-sugg-models} and the exact details on the benefits of
employing \eqref{eq:13} and \eqref{eq:14} will be discussed in
Section~\ref{sec:BregmanFramework}.\par
Let us also remark that the considered task in obtaining a representation
$AXA^{\top}$ from a given matrix $C$ bears strong similarities with matrix
factorisation problems, as described in \cite{C2007,XY2012,XYWZ2012a}, where
models, similar to ours, are also commonly used. Further, our models bear a
certain resemblance to \gls{BP} approaches, too. These usually seek solutions of
\begin{equation}
  \label{eq:15}
  \argmin_{x\in\Real^{n}}\set{\frac{1}{2}\norm{Bx-p}_{2}^{2} + \mu\norm{x}_{1}}
\end{equation}
for some given data $B\in\Real^{m,n}$, $p\in\Real^{m}$ and a positive
regularisation weight $\mu$. Such models are for example discussed in
\cite{HS2014,HS2015a}. Possible algorithms to obtain a solution to \eqref{eq:15}
include Lasso \cite{T1996b}, \gls{FBS} \cite{S2009a,S2010}, \gls{SB}
\cite{GO2009} and \gls{PDHG} \cite{CP2011,PC2011}. The \gls{PDHG} algorithm
requires an estimate of the eigenvalue spectrum of the matrix $B$ in order to
set the initialisation correctly. Depending on the size and the structure of the
matrix, such an estimate may be difficult to get. The \gls{FBS} can be very fast
if additional techniques to optimise the parameters of the framework are being
used. Split Bregman is very robust and does not need any special parameter
setups. This observation as well as its other convenient theoretical properties
make \gls{SB} our algorithm of choice to base our numerical strategies on.\par
% 
% ------------------------------------------------------------------------------
\subsection{Analysis of the Suggested Models}
\label{sec:analys-sugg-models}
% ------------------------------------------------------------------------------
% 
Let us analyse the suggested models and emphasise their properties in more
detail. The models proposed in \eqref{eq:10}-\eqref{eq:12} follow a hierarchy in
increasing complexity and accuracy on our model assumptions. All three models
are composed of two terms. The first term being a least squares formulation of
our assumption that \eqref{eq:7} should be fulfilled. The second term steers the
sparsity in the solution. The $\ell_{1}$ penaliser is a common choice to enforce
only few entries to be non-zero and easier to handle from an optimisation point
of view than directly penalising the number of non-zero entries. The latter
option would lead to a non-continuous combinatorial problem. The cost function
of each proposed model is convex, but not necessarily strictly convex.
Furthermore, only \eqref{eq:10} and \eqref{eq:12} possess cost functionals that
are always coercive.\par
The formulation in \eqref{eq:10} certainly resembles canonical $\ell_{1}$
regularised least-squares models most. The real valued non-negative scalar
$\lambda$ acts as a regularisation weight and influences the sparsity of the
solution. Larger parameter choices yield sparser solutions without allowing any
precise specification of the sparsity pattern. Its minimisation will be the
easiest to carry out but we expect the quality of the results to be inferior to
the other models, too. A diagonal matrix with only non-zero entries on its main
diagonal can already be regarded as sparse when compared to an arbitrary matrix.
Yet, we wish to have few non-zero entries along the main diagonal only. Thus the
model does not fully coincide with our assumptions.\par
Equation~\eqref{eq:11} allows a very fine grained tuning of the sparsity
structure. The weighting matrix $W$ with non-negative entries fulfils a similar
purpose as the parameter $\lambda$ in \eqref{eq:10} but allows a more
differentiated weighting of the individual entries. Setting all off-diagonal
entries to large values will favour solutions with diagonal structure. However,
a sparse structure along the main diagonal will still be difficult to achieve.
Suitable choices of the weighting matrix $W$ could also enable certain non
diagonal entries to be positive. Such set-ups would allow the handling of
correlated acoustic sources and could be the subject of future investigations.
The challenge lies in finding suitable weight distributions to accurately model
the correlation between individual sources. If all entries in $W$ are identical,
then \eqref{eq:11} coincides with \eqref{eq:10}.\par
Lastly, \eqref{eq:12} comes closest to our model assumptions that $X$ should be a
sparse diagonal matrix. The structural constraint is built explicitly into the
model. Such constrained optimisation tasks are usually more challenging than
their unconstrained counterparts, but in this case the cost considered in
\eqref{eq:12} comes closest to our model assumptions and is likely to yield the
best results. This model can be seen as the limiting case when all off-diagonal
elements of the weighting matrix $W$ in \eqref{eq:11} tend to $+\infty$.\par
Since \eqref{eq:10} represents a well studied model from the literature and since
\eqref{eq:11} and \eqref{eq:12} fit better to our model assumptions we will mostly
focus on the latter two formulations in this work.\par
% 
% ------------------------------------------------------------------------------
\section{Our Novel Solution Strategy}
\label{sec:SolvingStrategy}
% ------------------------------------------------------------------------------
% 
All suggested models are structured like
\begin{equation}
  \label{eq:16}
  \argmin_{x}\set{f(x) + g(Ax-b)}
\end{equation}
for some functions $f$, $g$, a matrix $A$, and a vector $b$. We note that this
setup also includes constrained optimisation tasks. It suffices to set
$g(x)=\iota_{S}(x)$, where $\iota_{S}$ is the indicator function of the set $S$
to force the solution to be inside the set $S$. Tasks like \eqref{eq:16} are
well suited to be solved by splitting schemes. Among the vast choice of existing
schemes we name the already mentioned \gls{FBS}, \gls{PDHG}, \gls{SB} as the
most popular ones. The \gls{FBS} algorithm probably belongs to the best studied
approaches and several extensions exist to further improve its efficiency. The
\gls{PDHG} approach excels in terms of speed, especially if it is used in
conjunction with preconditioning techniques such as those presented in
\cite{PC2011}. The \gls{SB} method is a very viable strategy with a thorough
convergence theory and sufficient flexibility to be applied to a large number of
distinct formulations.\par
In this work, we will adapt the \gls{SB} approach to our setting. Our choice is
motivated by the fact that the \gls{SB} method can be formulated such that it
applies to all our models and that it requires only a single parameter, which in
addition has a very intuitive interpretation.
%
% ------------------------------------------------------------------------------
\subsection{Mathematical Preliminaries}
\label{sec:Preliminaries}
% ------------------------------------------------------------------------------
%
In this section we briefly regroup a certain number of findings from the
literature that we will rely on in the forthcoming sections. These results stem
mostly from matrix (differential) calculus. We refer to
\cite{HJ1990,HJ1994,M2000a,MN2007,P1985,PP2008} for an in-depth discussion.\par
%
% ------------------------------------------------------------------------------
\subsubsection{Matrix Differential Theory}
\label{sec:matr-diff-theory}
% ------------------------------------------------------------------------------
%
The following results show that the space of complex valued matrices bears a
Hilbert space structure when equipped with the Frobenius norm. This observation
allows us to simplify certain expressions in the forthcoming sections and
provide a convenient framework to operate in.
\begin{definition}[Matrix Scalar Product \cite{PP2008}]
      Let $A = (a_{ij})$, $B = (b_{ij})$ be two $m\times n$ matrices over
      $\Complex$. We call \emph{matrix scalar product} the expression
      $\scprod{A}{B}$ defined by
      \begin{equation}
            \label{eq:17}
            \scprod{A}{B} \coloneqq \trace{A^{\top}B}
            = \sum_{i=1}^{m}\sum_{j=1}^{n} \overline{a_{ij}} b_{ij}\enspace{}.
      \end{equation}
\end{definition}
Here, the operator $\trace{\cdot}$ denotes the trace of a matrix. The matrix
scalar product is often called \emph{Frobenius scalar product} because it
induces the Frobenius norm.
\begin{lemma}
      Let $A = (a_{ij})$ be a $m\times n$ matrix over $\Complex$. The Frobenius
      norm is related to the matrix scalar product in the following way:
      \begin{equation}
            \label{eq:18}
            \norm{A}_{\Fro}^{2}
            \coloneqq \sum_{i=1}^{m}\sum_{j=1}^{n} \abs{a_{ij}}^{2}
            = \sum_{i=1}^{m}\sum_{j=1}^{n} \overline{a_{ij}} a_{ij}
            = \scprod{A}{A}
      \end{equation}
\end{lemma}
\begin{proposition}
      The matrix scalar product is a scalar product in the proper mathematical
      sense, i.e.\ all properties that are known to hold for general scalar
      products also apply to the matrix scalar product. In particular, we have
      the following rules for complex valued matrices $A$, $B$, and $C$:
      \begin{gather}
            \label{eq:19}
            \norm{A\pm{}B}_{\Fro}^{2} =
            \norm{A}_{\Fro}^{2} \pm{}2\scprod{A}{B} + \norm{B}_{\Fro}^{2}\\
            \label{eq:20}
            \scprod{A\, B}{C} = \scprod{B}{A^{\top}C}
            \quad{}\text{and}\quad
            \scprod{A}{B\, C} = \scprod{B^{\top}A}{C}
      \end{gather}
      \begin{proof}
            The proofs are straightforward. Equation~\eqref{eq:19} follows from
            the bilinearity of the scalar product. In order to show \eqref{eq:20}
            it is helpful to use the following result from \cite{PP2008}
            (Equation~(16)):
            \begin{equation}
                  \label{eq:21}
                  \trace{A\, B\, C} = \trace{B\, C\, A} = \trace{C\, A\, B}
            \end{equation}
            \qed
      \end{proof}
\end{proposition}
Splitting schemes such as the \gls{SB} method often assume that the cost
function to be minimised can be decomposed into simpler (convex) terms which
will be optimised in certain ways alternatingly. For performance reasons it is
beneficial if these optima have closed form solutions. Such representations can
usually be derived from first order optimality conditions. The next propositions
state a certain number of findings related to matrix valued derivatives that
will prove to be helpful. A thorough discussion on the notion of matrix
differentials and their derivatives can be found in \cite{MN2007}. As we mainly
make use of the results documented in the aforementioned references we only
state the necessary formulas here.
\begin{corollary}
      We have $\scprod{A\,X\, B^{\top}}{Y}=\scprod{X}{A^{\top} Y\, B}$.
      \begin{proof}
            Follows from \eqref{eq:21}.
      \end{proof}
\end{corollary}
Computing derivatives with respect to matrices requires to take the structure of
the matrix into account. In a symmetric matrix off-diagonal entries appear
always twice for example. We denote the derivative with respect to a matrix $X$
with $\frac{\mathrm{d}}{\mathrm{d}X}$. If $X$ is assumed to have a specific
structure we will mention it explicitly.
\begin{proposition}
      \label{thm:2}
      We have the following identities for an arbitrary unstructured matrix~$X$:
      \begin{align}
        \label{eq:22}
        \frac{\mathrm{d}}{\mathrm{d} X}\trace{(A\, X\, B + C)^{\top}(A\, X\, B + C)}
        &= \frac{\mathrm{d}}{\mathrm{d} X} \norm{A\, X\, B + C}_{\Fro}^{2} \\
        \label{eq:23}
        &= 2\, A^{\top}(A\, X\, B+C)B^{\top} \\
        \label{eq:24}
        \frac{\mathrm{d}}{\mathrm{d} X}\norm{W\circ{}X-C}_{\Fro}^{2}
        &= 2\, W\circ\left( W\circ{}X-C \right)
      \end{align}
      From \eqref{eq:22} and \eqref{eq:23} we may conclude that necessary
      optimality conditions for a minimiser of $\frac{1}{2}\norm{A\, X\,
        A+C}_{\Fro}^{2}$ are given by $A^{\top}(A\, X\, A+C)A^{\top}=0$.
      \begin{proof}
            Eq.~\eqref{eq:23} is stated explicitly in \cite{PP2008} as
            Eq.~(108), whereas \eqref{eq:24} follows from the chain rule
            (Theorem 12, Chapter 5 in \cite{MN2007}) and the differential of the
            Hadamard product (Eq.~(17), Chapter 8 in \cite{MN2007}).
      \end{proof}
\end{proposition}
If we require that the matrix is structured, then other identities hold
\cite{M2000a}. When working with symmetric (resp.\ diagonal) matrices it seems
intuitive to expect derivatives to yield symmetric (resp.\ diagonal) matrices as
well. Especially, if we restrict our attention to diagonal matrices, then we
obtain the following equality.
\begin{proposition}
      \label{thm:3}
      We have the following equality for a \emph{diagonal} matrix $X$:
      \begin{equation}
            \label{eq:25}
            \frac{\mathrm{d}}{\mathrm{d}X}
            \frac{1}{2}\norm{A\, X\, B-C}_{\Fro}^{2} =
            \left( A^{\top}\left( A\, X\, B-C \right)B^{\top} \right)\circ{}I
      \end{equation}
      where $I$ is the identity matrix.
      \begin{proof}
            Follows from the formulas established in Proposition~\ref{thm:2}
            under consideration of the restrictions stated in \cite{PP2008}
            (Section 2.8).
      \end{proof}
\end{proposition}
By comparing \eqref{eq:22} and \eqref{eq:25} we clearly see the importance of
taking the structure of the underlying matrices into account. The derivative in
\eqref{eq:25} only coincides with the derivative in \eqref{eq:22} along the main
diagonal. A similar finding can be deduced for symmetric matrices, too.\par
%
% ------------------------------------------------------------------------------
\subsubsection{Matrix Valued Proximal Operators}
\label{sec:matr-valu-prox}
% ------------------------------------------------------------------------------
%
As already mentioned, splitting algorithms decompose the original task into
simpler subproblems. Some of these subtasks coincide with proximal operations
which go back to Moreau \cite{M1965}. We also refer to \cite{GK2002} for a
detailed analysis of their properties. The proximal operator of the $\ell_{1}$
norm will play a significant role in our framework. It is a well known operator
for signal processing tasks, where it is often referred to as soft shrinkage.
\begin{definition}[Soft Shrinkage]
      The soft shrinkage operator $\shrink_{\lambda}(b)$ with parameter
      $\lambda>0$ solves the optimisation problem
      \begin{equation}
            \label{eq:26}
            \argmin_{x\in\Real^{n}}
            \set{\lambda\norm{x}_{1} + \frac{1}{2}\norm{x-b}_{2}^{2}}
      \end{equation}
\end{definition}
For real valued data the soft shrinkage operator has a closed form
representation. We provide this result with a detailed proof since we could not
find a source with a satisfying derivation.
\begin{lemma}
      \label{thm:4}
      Let $b\in\Real^{n}$. The solutions of \eqref{eq:26} are given by
      \begin{equation}
            \label{eq:27}
            \begin{split}
                  \left( \shrink_{\lambda}(b) \right)_{i}
                  &\coloneqq \text{sgn}\left(b_{i}\right)
                  \max\left(\abs{b_{i}}-\lambda, 0\right) \\
                  &=
                  \begin{cases}
                        b_{i} + \lambda, &b_{i}<-\lambda \\
                        0, &b_{i}\in\left[-\lambda,\lambda\right] \\
                        b_{i} - \lambda, &b_{i}>\lambda
                  \end{cases}
            \end{split}
      \end{equation}
      \begin{proof}
            First, we note that \eqref{eq:26} decouples into $n$ independent
            optimisation problems of the form
            \begin{equation}
                  \label{eq:28}
                  \argmin_{x_{i}\in\Real}
                  \set{\lambda\abs{x_{i}} + \frac{1}{2}(x_{i}-b_{i})^{2}}
                  \quad i=1, \ldots, n
            \end{equation}
            Thus, it suffices to solve the simpler one dimensional optimisation
            task and to rearrange the result componentwise. We note that the energy in
            \eqref{eq:28} is a strictly convex function and thus every local minimum must
            necessarily also be a global minimum. Minimisers $\overline{x}_{i}$ must fulfil
            the first order optimality condition $0\in
            \lambda \partial\left(\abs{\cdot}\right) \left( \overline{x}_{i} \right) +
            \left(\overline{x}_{i}-b_{i}\right)$, where $\partial\left(\abs{\cdot}\right)
            \left( \overline{x}_{i} \right)$ is the subdifferential of the absolute value
            function evaluated at $\overline{x}_{i}$. Let $q
            \in \partial\left(\abs{\cdot}\right) \left( \overline{x}_{i} \right)$ be any
            subgradient, then we obtain
            \begin{equation}
                  0 = \lambda q + \overline{x}_{i} - b\
                  \Leftrightarrow\
                  q = \lambda^{-1} b - \lambda^{-1} \overline{x}_{i}
            \end{equation}
            There are three possible choices for $q$. If $\overline{x}_{i} > 0 $
            then $q=1$, if $\overline{x}_{i}<0$ then $q=-1$ and if $\overline{x}_{i}=0$,
            then $q\in\left[-1,1\right]$. Now, if $\lambda^{-1} b > 1$ then
            $\overline{x}_{i}$ must be positive because $q$ cannot be larger than 1. But for
            positive $\overline{x}_{i}$, $\partial\left(\abs{\cdot}\right)(\overline{x}_{i})
            = \set{q} = \set{1}$ and $\overline{x}_{i}$ must be $b-\lambda$. On the other
            side, if $\lambda^{-1} b < -1$, then $\overline{x}_{i}$ must be negative, since
            $q$ cannot be smaller than $-1$. For negative $\overline{x}_{i}$,
            $\partial\left(\abs{\cdot}\right)(\overline{x}_{i}) = \set{q} = \set{-1}$ and
            $\overline{x}_{i}$ is given by $b+\lambda$. Now assume $\lambda^{-1} b \in
            \left[0,1\right]$, then $-\lambda^{-1} \overline{x}_{i} = q - \lambda^{-1} b$
            must hold. If we suppose $q > \lambda^{-1} b$, then $\overline{x}_{i}$ would
            have to be strictly negative. However, for strictly negative $\overline{x}_{i}$,
            $q=-1$ and we have a contradiction. In the same way we cannot have $q <
            \lambda^{-1} b$ since then we must have $\overline{x}_{i}>0$ and $q=1$. It
            follows that $q=\lambda^{-1} b$ is the only remaining choice. This implies
            $\overline{x}_{i}=0$. The case $\lambda^{-1} b \in \left[-1,0\right]$ leads by
            identical reasoning to the same result $\overline{x}_{i}=0$ and $q= \lambda^{-1}
            b$. Regrouping all solutions completes the proof.
            \qed
      \end{proof}
\end{lemma}
Lemma~\ref{thm:4} yields a closed form description of the proximal mapping for
real valued arguments only. However, our task at hand has complex valued data.
Ignoring the imaginary parts in the computations would lead to significant
errors. The following theorem shows how to extend the previous finding. We note,
that it relies on our choice for the complex valued norms from \eqref{eq:13} and
\eqref{eq:14}. For other choices it might be difficult or even not possible at
all to get such a closed form representation.
\begin{theorem}
      \label{thm:5}
      Let $B\in\Complex^{n,m}$. Then, we have
      \begin{equation}
            \label{eq:29}
            \argmin_{X\in\Complex^{n,m}} \set{
              \lambda\norm{X}_{1} + \frac{1}{2}\norm{X-B}_{\Fro{}}^{2}
            } =
            \left(
                  \shrink_{\lambda}\left(  \Rpart (B) \right) +
                  \imath\cdot \shrink_{\lambda}\left(  \Ipart (B) \right)
            \right)_{ij}
      \end{equation}
      with $\imath$ denoting the complex unit.
      \begin{proof}
            Due to our convention in \eqref{eq:13} and \eqref{eq:14}, the
            minimisation in \eqref{eq:29} not only decouples into $nm$ smaller
            optimisation problems, each in a single scalar valued complex
            unknown, but we can further optimise real and imaginary part
            independently:
            \begin{multline}
                  \label{eq:30}
                  \argmin_{X\in\Complex^{n,m}} \set{
                    \lambda\norm{X}_{1} + \frac{1}{2}\norm{X-B}_{\Fro{}}^{2}
                  } = \\
                  \argmin_{X\in\Complex^{n,m}} \set{\lambda
                  \sum_{i,j} \abs{\Rpart (X_{ij})} + \abs{\Ipart (X_{ij})}
                  +
                  \frac{1}{2}\left( \Rpart (X_{ij}-B_{ij})^{2} +
                        \Ipart (X_{ij}-B_{ij})^{2} \right)}
            \end{multline}
            Thus, we obtain $2mn$ optimisation problems that can be solved
            independently by making use of Lemma~\ref{thm:4}.
            \qed
      \end{proof}
\end{theorem}
We note that the previous theorem naturally extends the case of real valued
matrices as well as the setting with vector valued data. It also represents an
important finding for our numerical scheme. The componentwise handling of the
data allows an efficient scaling of our numerical solvers to large parallel
processing facilities such as \glspl{GPGPU}.\par
% 
% ------------------------------------------------------------------------------
\subsection{The Bregman Framework}
\label{sec:BregmanFramework}
% ------------------------------------------------------------------------------
% 
We present a short review of the (split) Bregman algorithm, which will present
the basis of our forthcoming developments. There exists a large family of
algorithms that could be considered to form the Bregman framework. The most
prominent representatives would be the \emph{standard Bregman iteration}, as
developed by Osher and colleagues\ \cite{Osher2005}, the \emph{linearised
Bregman algorithm} by Cai et al.\ \cite{Cai2009}, and the \emph{split Bregman
method} by Goldstein and Osher \cite{GO2009}. Bregman himself \cite{Bregman1967}
originally wanted to describe non-orthogonal projections onto convex sets and
derived, as by-product, an iterative scheme to minimise certain smooth and
convex functions under linear constraints. We refer to \cite{BB1997} for a
deeper analysis of these projections. Let us also remark that the Bregman
algorithms are related to many other popular strategies from the literature and
a certain number of equivalences with proximal methods have been discovered in
the past. An extensive discussion on these relationships can for example be
found in \cite{S2010}. Convergence considerations for various setups can, in
addition to the already cited works, also be found in
\cite{Burger2006,BRH07,Cai2009a,GO2009,BSB2011}. Finally, various strategies to
improve the computational performance of the split Bregman scheme has been
discussed in \cite{GTBS2015, GDSB2014a, GTSJ2014a, KCSB2015} and the references
therein.\par
% 
% ------------------------------------------------------------------------------
\subsubsection{The Split Bregman Method}
\label{sec:SplitBregmanFramework}
% ------------------------------------------------------------------------------
% 
The \gls{SB} method has originally been developed to solve convex $\ell_{1}$
penalised optimisation tasks. In its simplest form it aims to find a solution of
equations such as \eqref{eq:16} with convex but not necessarily smooth
functions $f$ and $g$. The underlying idea is to introduce a slack variable $z$
and to reformulate the original task in the equivalent form
\begin{equation}
  \label{eq:31}
  \begin{gathered}
    \argmin_{x,\ z}\set{f(x) + g(z)}\\
    \text{such that}\quad z = Ax-b
  \end{gathered}
\end{equation}
and to carry out the minimisations with resp.\ to $x$ and $z$ in an alternating
manner. If there exist closed form representations of the proximal mappings
\begin{align}
  \label{eq:32}
  \argmin_{x} &\set{f(x) + \frac{\mu}{2}\norm{A x-p_{1}}_{2}^{2}} \\
  \label{eq:33}
  \argmin_{z} &\set{g(z) + \frac{\mu}{2}\norm{z-p_{2}}_{2}^{2}}
\end{align}
for arbitrary vectors $p_{1}$, $p_{2}$, and arbitrary positive $\mu$, then the
\gls{SB} algorithm can be carried out very fast and avoid possible pitfalls that
would occur in the formulation of \eqref{eq:16} due to the potential
non-differentiability. The complete formulation of the \gls{SB} algorithm to
solve \eqref{eq:31} is given in Algorithm~\ref{alg:SplitBregmanGeneric}.
\begin{algorithm}
  \caption{Split Bregman Algorithm for solving \eqref{eq:31}}
  \label{alg:SplitBregmanGeneric}
  \DontPrintSemicolon
  Initialise $x^{[0]}$ such that $0\in \partial \left( g \right)(x^{[0]})$, $z^{[0]} = 0$,
  $\hat{b} = b$ and $\mu > 0$ arbitrarily\;
  \Repeat{convergence}{
    \Repeat{convergence}{
      $\displaystyle{}x^{[k+1]} =
      \argmin_{x}\set{ f(x) + \frac{\mu}{2}\norm{Ax - z^{[k]} - \hat{b}}_{2}^{2}}$
      \label{alg:SBItMin1}\;
      $\displaystyle{}z^{[k+1]} =
      \argmin_{x}\set{ g(z) + \frac{\mu}{2}\norm{Ax^{[k+1]} - z - \hat{b}}_{2}^{2}}$
      \label{alg:SBItMin2}\;
    }
    $\hat{b} = \hat{b} - \left( Ax^{[k+1]} - b - z^{[k+1]} \right)$\;
  }
\end{algorithm}
We remark that the initialisation stated in
Algorithm~\ref{alg:SplitBregmanGeneric} asserts convergence for any choice $\mu
> 0$, although the convergence speed may be affected. In practical setups it may
be difficult to find an initial value for $x$ such that $0\in \partial g(x)$. An
asset of the \gls{SB} scheme is that it is still likely to converge with
arbitrary initialisations. The two minimising steps in Lines~\ref{alg:SBItMin1}
and \ref{alg:SBItMin2} represent a minimising strategy to solve
\begin{equation*}
      \argmin_{x,\ z}\set{f(x) + g(z) + \frac{\mu}{2}\norm{Ax - z - \hat{b}}_{2}^{2}}
\end{equation*}
The original authors of the \gls{SB} algorithm suggest a single sweep in
\cite{GO2009}. Our experience suggests that a few more alternating minimisations
may be worthwhile for certain applications.\par
In the forthcoming paragraphs we will show how to fit the formulations from
Algorithm~\ref{alg:SplitBregmanGeneric} to our models. The main difficulty is
the fact that our optimisation is being carried out with respect to complex
valued (diagonal) matrices.
%
% ------------------------------------------------------------------------------
\subsection{Optimising Without Structural Constraints}
\label{sec:OptSplitBregWeighted}
% ------------------------------------------------------------------------------
% 
In this section we derive an iterative strategy to solve \eqref{eq:10} and
\eqref{eq:11}. Both models only differ in a weighting matrix and can be handled
in the same algorithmic way. We suggest to apply the following mapping
\begin{align}
  \label{eq:34}
  f(X) &\coloneqq \frac{1}{2}\norm{AXA^{\top}-C}_{\Fro}^{2} \\
  \label{eq:35}
  g(X) &\coloneqq \lambda \norm{W \circ X}_{1}
\end{align}
where, depending on whether we consider \eqref{eq:10} or \eqref{eq:11}
(resp.~\eqref{eq:12} in the forthcoming section), either $\lambda=1$, or $W$ is
the matrix having 1 as entry in each position. In addition, we introduce the
slack variable by requiring that $X=D$. A straightforward application of the
\gls{SB} algorithm (as presented in \eqref{eq:31}) leads now to the following
iterative formulation
\begin{equation}
  \label{eq:36}
  \begin{split}
    (X,D)^{[k+1]} &= \argmin_{X,\; D\in\Real^{m,m}} \set{ 
      f(X) + g(D) + \frac{\mu}{2} \norm{D - X - B^{[k]}}_{\Fro}^{2}
    } \\
    B^{[k+1]} &= B^{[k]} - D^{[k+1]} + X^{[k+1]}
  \end{split}
\end{equation}
The minimisation with respect to $(X, D)$ is done in an alternating manner. The
optimisation with respect to $X$ reduces to solving the following least squares
problem
\begin{equation}
  \label{eq:37}
  X^{[k+1]} = \argmin_{X\in\Real^{m,m}}\set{
    \frac{1}{2}\norm{AXA^{\top} - C}_{\Fro}^{2} + 
    \frac{\mu}{2}\norm{X - \left( D^{[k]} - B^{[k]} \right)}_{\Fro}^{2}
  }
\end{equation}
while the minimisation of $D$ decouples into $m^{2}$ optimisations. Indeed, we
have to solve
\begin{equation}
  \label{eq:38}
  D^{[k+1]} = \argmin_{D\in\Real^{m,m}}\set{
    \lambda\norm{W\circ D}_{1} + 
    \frac{\mu}{2}\norm{D - X - B^{[k]}}_{\Fro}^{2}
  }
\end{equation}
which can be reduced to component-wise soft shrinkage operations as shown in
\eqref{eq:29}. The necessary optimality conditions for \eqref{eq:37} are given
by the linear system
\begin{equation}
  \label{eq:39}
  A^{\top} \left( AXA^{\top} - C \right) A + \mu (X - D^{[k]} - B^{[k]}) = 0
\end{equation}
The minimisation in \eqref{eq:38} can be expressed in terms of the soft
shrinkage operator and a closed form expression for $D_{ij}$ is given by
\begin{equation}
  \label{eq:40}
  \begin{multlined}[c][0.7\linewidth]
        D_{ij}^{[k+1]} = \shrink_{\frac{\lambda W_{ij}}{\mu}} \left( \Rpart
              \left( X_{ij}^{[k+1]} + B_{ij}^{[k]} \right) \right) + \\
        \imath \cdot \shrink_{\frac{\lambda W_{ij}}{\mu}} \left( \Ipart \left(
                    X_{ij}^{[k+1]} + B_{ij}^{[k]} \right) \right)
  \end{multlined}
\end{equation}
We emphasise the importance that $W$ is a real-valued matrix with non-negative
entries only. Otherwise it would not be possible to extract its entries and to
combine them with the parameter $\mu$ in the shrinkage formula.\par
While the minimisation with respect to $D$ can be carried out very efficiently,
this is not necessarily the case for the minimisation with respect to $X$.
Indeed, the occurring linear system has a very large and dense matrix. Solving
this system with standard methods from the literature would lead to a strategies
that are prohibitive both in terms of memory and run time. Instead, we opt to
solve \eqref{eq:37} through a simple gradient descent scheme. A single descent
step can easily be derived and is given by
\begin{equation}
  \label{eq:41}
  X^{[k,j+1]} = X^{[k,j]} - \alpha \left(
    A^{\top}\left( AX^{[k,j]}A^{\top}-C \right) A + 
    \mu \left( X^{[k,j]}-D^{[k]}+B^{[k]} \right)
  \right)
\end{equation}
where $\alpha$ is a step size that must be chosen sufficiently small to assert a
decrease in energy. We note that the best value for $\alpha$ may be obtained
numerically through a simple 1D optimisation and that the best value can also be
estimated analytically. We refer to \cite{NW2006} for an overview of several
efficient numerical strategies and to Proposition~\ref{thm:OptStepSize} for the
computation of the optimal step size. The complete algorithm is given in
Algorithm~\ref{alg:SplitBregmanNoStructure}. Let us also remark that the two
nested inner loops in the algorithm do not need to be carried out until full
convergence. For many applications already a single iteration step is usually
enough: In \cite{GO2009} the authors suggested to use a single Gauß-Seidel
iteration to solve the occurring linear system. The influence of inaccuracies in
the Bregman iteration has also been studied in \cite{YO2013}. In general the
\gls{SB} scheme is very robust against inaccuracies. The benefit of reducing the
number of iterations will usually yield significant increase in performance.
\begin{algorithm}[bt]
  \DontPrintSemicolon
  \caption{Split Bregman for solving \eqref{eq:10} or \eqref{eq:11}}
  \label{alg:SplitBregmanNoStructure}
  \KwIn{$A$, $C$, $\lambda$, $W$, $\alpha$}
  \KwOut{Solution of \eqref{eq:10} resp.\ \eqref{eq:11}.}
  Initialise $X=0$, $D=0$ and $B = 0$\;
  % ----------------------------------------------------------------------------
  \Repeat{convergence of $X$, $D$ and $B$}{
    \label{alg:SBMloop1}
    Set $\hat{X} = X$ and $\hat{D} = D$\;
    % --------------------------------------------------------------------------
    \Repeat{convergence towards $\hat{X}^{*}$ and $\hat{D}^{*}$}{
      \label{alg:SBMloop2}
      Set $\overline{X} = \hat{X}$\;
      % ------------------------------------------------------------------------
      \Repeat{convergence towards $\overline{X}^{*}$}{\label{alg:SBMloop3}
        Compute $\alpha$ according to Proposition~\ref{thm:OptStepSize}.\;
        $\overline{X} = \overline{X} - 
        \alpha \left( A^{\top} \left( A \overline{X} A^{\top} - C \right) A +
              \lambda \left( \overline{X} - \hat{D} + B \right)
        \right)$\;
      }
      % ------------------------------------------------------------------------
      $\hat{D} = \shrink_{\frac{\mu W}{\lambda}}
      \left( \overline{X}^{*} + B \right)$\;
    } % ------------------------------------------------------------------------
    Set $X = \hat{X}^{*}$, $D = \hat{D}^{*}$ and $B = B - D + X$\;
  }
\end{algorithm}
The first two nested loops in Algorithm~\ref{alg:SplitBregmanNoStructure}
correspond in a one-to-one manner to the loops in
Algorithm~\ref{alg:SplitBregmanGeneric}. The third loop beginning in
Line~\ref{alg:SBMloop3} represents the gradient descent scheme explained in
\eqref{eq:41}. It corresponds to the minimisation with respect to $x$ in
Line~\ref{alg:SBItMin1} in Algorithm~\ref{alg:SplitBregmanGeneric}.\par
The following proposition states the optimal value for the gradient descent
scheme, in the sense that in each iteration the decrease in the energy is
maximal.
\begin{proposition}
      \label{thm:OptStepSize}
      The optimal value for the step size $\alpha$ for the gradient descent
      scheme from \eqref{eq:41} is given by:
      \begin{equation}
            \label{eq:42}
            \alpha =
            \frac{\norm{G(X^{[k,j]})}_{\Fro}^{2}}{\norm{AG(X^{[k,j]})A^{\top}}_{\Fro}^{2}
              + \mu\norm{G(X^{[k,j]})}_{\Fro}^{2}}
      \end{equation}
      where $G$ is the first order derivative with respect to $X$ of the cost functional 
      \begin{equation}
            \label{eq:43}
            E(X) \coloneqq \frac{1}{2}\norm{AXA^{\top}-C}_{\Fro}^{2} +
            \frac{\mu}{2}\norm{D^{[k]}-X-B^{[k]}}_{\Fro}^{2}
      \end{equation}
      \begin{proof}
            Let $G$ be the gradient of our energy $E$, i.e.\
            \begin{equation}
                  \label{eq:44}
                  G(X) = A^{\top} (A\; X\; A^{\top} - C ) A + \mu ( X - D^{[k]} + B^{[k]})
            \end{equation}
            In order to maximise the energy decrease in a single gradient
            descent step we wish to solve
            \begin{equation}
                  \label{eq:45}
                  \argmin_{\alpha}\set{E\left( X^{[k,j]} - \alpha\; G(X^{[k,j]}) \right)}
            \end{equation}
            Using the definition of $E$, we are led to the following expression,
            that needs to be minimised with respect to $\alpha$:
            \begin{equation}
                  \label{eq:46}
                  \frac{1}{2}\norm{ \alpha 
                    AG(X^{[k,j]})A^{\top}-( AXA^{\top}-C )}_{\Fro}^{2} +
                  \frac{\mu}{2}\norm{\alpha
                    G(X^{[k,j]})-( X^{[k,j]}-D+B^{(k)} )}_{\Fro}^{2}
            \end{equation}
            Using the properties of the matrix scalar product, it follows that the latter
            expression can be rewritten as
            \begin{equation}
                  \label{eq:47}
                  \begin{multlined}[0.8\textwidth]
                        E(X^{[k,j]}) + \alpha^{2}
                        \left( 
                              \frac{1}{2} \norm{AG(X^{[k,j]})A^{\top}}_{\Fro}^{2} +
                              \frac{\mu}{2} \norm{G(X^{[k,j]})}_{\Fro}^{2} \right) \\
                        - \alpha \scprod{AG(X^{[k,j]})A^{\top}}{AX^{[k,j]}A^{\top}-C}
                        - \alpha \mu \scprod{G(X^{[k,j]})}{X^{[k,j]}-D+B^{(k)}}
                  \end{multlined}
            \end{equation}
            Thus, the optimal $\alpha$ is given by
            \begin{equation}
                  \label{eq:48}
                  \begin{split}
                        \alpha &= \frac{
                          \scprod{AG(X^{[k,j]})A^{\top}}{AX^{[k,j]}A^{\top}-C} + 
                          \mu \scprod{G(X^{[k,j]})}{X^{[k,j]}-D+B^{(k)}}}
                        {
                          \norm{AG(X^{[k,j]})A^{\top}}_{\Fro}^{2} +
                          \mu\norm{G(X^{[k,j]})}_{\Fro}^{2}
                        } \\
                        &= \frac{
                          \scprod{G(X^{[k,j]})}{A^{\top}AX^{[k,j]}A^{\top}A-A^{\top}CA +
                            \mu (X^{[k,j]}-D+B^{(k)})}}
                        {
                          \norm{AG(X^{[k,j]})A^{\top}}_{\Fro}^{2} +
                          \mu\norm{G(X^{[k,j]})}_{\Fro}^{2}
                        } \\
                        &= \frac{
                          \norm{G(X^{[k,j]})}_{\Fro}^{2}
                        }{
                          \norm{AG(X^{[k,j]})A^{\top}}_{\Fro}^{2} +
                          \mu\norm{G(X^{[k,j]})}_{\Fro}^{2}
                        }
                  \end{split}
            \end{equation}
            \qed
      \end{proof}
\end{proposition}
% 
% ------------------------------------------------------------------------------
\subsection{Structured Split Bregman Model}
\label{sec:OptSplitBregStruct}
% ------------------------------------------------------------------------------
% 
We now consider our third model from \eqref{eq:12}. Our approach is similar to
the approach presented in Section~\ref{sec:OptSplitBregWeighted}. This time we
additionally enforce a diagonal structure of $X$ by restricting the
optimisation with respect to $X$ to matrices having the desired structure. Our
functional takes the form
\begin{equation}
  \label{eq:49}
  \argmin_{X}\set{
    \frac{1}{2}\norm{AXA^{\top}-C}_{\Fro}^{2} + \mu\norm{X}_{1}}
\end{equation}
where $\mu$ is again a sparsity inducing regularisation weight. We abstain from
using a more fine-grained weighting matrix $W$ since we do not have any prior
information on the sparsity pattern. Nevertheless, we remark that such a
modification would be possible.\par
This time the minimisation is not done over the set of all matrices but merely
over the set of diagonal matrices. Our consideration is that the reduced search
space improves performance and yields optima closer to our expectations.\par
The minimisation is done with the help of the \gls{SB} approach. Similarly as
before, we introduce a dummy variable $D$ and obtain
\begin{equation}
  \label{eq:50}
  \begin{split}
    (X,D)^{[k+1]}
    &=\argmin_{X,D}\set{
      \frac{1}{2} \norm{AXA^{\top}-C}_{\Fro}^{2} +
      \mu\norm{D}_{1}+\frac{\lambda}{2} \norm{D-X-B^{[k]}}_{\Fro}^{2}} \\
    B^{[k+1]} &= B^{[k]} - D^{[k+1]} + X^{[k+1]}
  \end{split}
\end{equation}
The minimisation with respect to $X$ and $D$ is done in an alternating manner.
The optimisation with respect to $D$ is identical to the previous section where
we set $W_{ij}=1$ for all $i$ and $j$. If $D$ is initialised as a diagonal
matrix, then no special consideration need to be done to preserve this structure
provided that $X$ is initialised as a diagonal matrix, too. The minimisation
with respect to $X$ cannot be reduced to a linear system anymore since we have
introduced a structure preserving constraint. Nevertheless, an explicit form of
the gradient of the energy to be minimised is available. Using the findings from
Proposition~\ref{thm:3} we conclude that the gradient of our cost function in
\eqref{eq:50} with respect to the \emph{diagonal matrix} $X$ can be expressed as
\begin{equation}
  \label{eq:51}
  \left( A^{\top} \left( AXA^{\top}-C \right) A \right) \circ I +
  \lambda \left( X-D^{(k)}+B^{(k)} \right) \circ I
\end{equation}
The two Hadamard products with the identity matrix in the previous formula stem
from the fact that we compute the derivative with respect to diagonal matrices.
They ensure that the diagonal structure of the matrix $X$ is preserved
throughout the whole iterative process.\par
The full algorithm is given in Algorithm~\ref{alg:SplitBregmanStructure}. Also
here it is possible to optimise the step size for the gradient descent scheme in
a similar manner as in Proposition~\ref{thm:OptStepSize}. 
\begin{algorithm}[tb]
  \DontPrintSemicolon
  \caption{Split Bregman for solving \eqref{eq:12}}
  \label{alg:SplitBregmanStructure}
  \KwIn{$A$, $C$, $\alpha$, $\mu$}
  \KwOut{Solution of \eqref{eq:12}.}
  Initialise $X=0$, $D=0$ and $B=0$\;
  % ----------------------------------------------------------------------------
  \Repeat{convergence of $X$, $D$ and $B$}{
    \label{alg:SSBMloop1-2}
    Set $\hat{X} = X$ and $\hat{D} = D$\;
    % --------------------------------------------------------------------------
    \Repeat{convergence towards $\hat{X}^{*}$ and $\hat{D}^{*}$}{
      \label{alg:SSBMloop2-2}
      Set $\overline{X} = \hat{X}$\;
      % ------------------------------------------------------------------------
      \Repeat{convergence towards $\overline{X}^{*}$}{
        \label{alg:SBMloop3-2}
        Compute optimal $\alpha$.\;
        $\overline{X} = \overline{X} -
        \alpha \left( A^{\top} \left(  A \overline{X} A^{\top}-C \right) A +
              \lambda \left(\overline{X} - \hat{D} + B \right)
        \right) \circ I$\;
      }
      % ------------------------------------------------------------------------
      $\hat{D} = \shrink_{\frac{\mu}{\lambda}}
      \left( \hat{X} + B \right)$\;
    }
    % --------------------------------------------------------------------------
    Set $X = \hat{X}^{*}$, $D = \hat{D}^{*}$ and $B = B - D + X$\;
  }
\end{algorithm}
%
% ------------------------------------------------------------------------------
\subsection{Post Processing Steps}
\label{sec:post-proc-steps}
% ------------------------------------------------------------------------------
%
Numerical experiments show that sometimes it is necessary to post-process the
obtained signals to remedy certain deficiencies that stem from the experimental
nature of the physical setup. The following simple strategy has proven to be
efficient and reliable. We formulate it for the results obtained from
Section~\ref{sec:OptSplitBregStruct}, but an analogous application to the
results obtained with the algorithm from Section~\ref{sec:OptSplitBregWeighted}
are clearly possible, too.
\begin{enumerate}
\item Remap the diagonal entries from the matrix $X$ to their actual positions
      in 2D (resp.\ 3D) space. Thus, we obtain tuples $(i, j, x_{ij})$.
\item Apply a k-means clustering \cite{M1967} approach to partition the tuples
      into distinct sets. If the exact number of clusters is unknown then it can
      either be estimated with the strategy from \cite{S2010a} or by using other
      popular strategies commonly used in the clustering context, such as
      silhouette coefficients \cite{R1987} or GAP statistics \cite{TWH2001}.
\item Use the centroid position of each cluster as source position.
\item Sum up all source strengths from a cluster to obtain the corresponding
      source strength.
\end{enumerate}
Our experiments have shown that the sum of all entries of our computed solutions
always coincided with the sum of all entries of the ground truth solution. This
observation motivates our suggestion to sum all entries of a cluster. Further
knowledge on the solution could also be introduced into the post processing. As
such, it would for example be possible to introduce a minimal distance between
two peaks. This could help in handling noise, as it tends to cluster small
fluctuations around a single large peak.
%
% ------------------------------------------------------------------------------
\section{Numerical Evaluation}
\label{sec:NumericalEvaluation}
% ------------------------------------------------------------------------------
%
We evaluate our models and show the benefits of the additional post processing.
For each tested model we also discuss the influence of noisy data and
inaccuracies in the measurements. The findings are depicted in their
corresponding sections.\par
\begin{figure}
      \centering
      \includegraphics{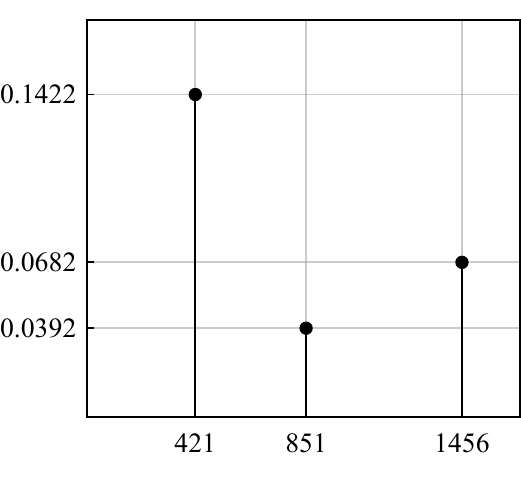}
      \hspace*{2em}
      \raisebox{0.05\height}{
        \begin{tikzpicture}[scale=9.25,>=latex]
              \def\xym{0.22}
              \draw[step=0.1, gray, very thin] (-0.205,-0.205) grid (0.205,0.205);
              \draw[black,only marks, mark=o,mark options={scale=0.2}] plot (-0.1,-0.1);
              \draw[black,only marks, mark=o,mark options={scale=0.2}] plot (0.15, 0.0);
              \draw[black,only marks, mark=o,mark options={scale=0.2}] plot ( 0.0, 0.1);
              \draw[black, fill=white] (-0.1,-0.1) node[anchor=north east] {$0.1422$};
              \draw[black, fill=white] (0.15, 0.0) node[anchor=north, yshift=-1mm] {$0.0682$};
              \draw[black, fill=white] ( 0.0, 0.1) node[anchor=north west] {$0.0392$};
              \foreach \x in {-0.2,-0.1,0,0.1,0.2}
              \draw (\x,-\xym+0.01) -- +(0,-0.01) node [below] {\x};
              \foreach \y in {-0.2,-0.1,0,0.1,0.2}
              \draw (-\xym+0.01,\y) -- +(-0.01,0) node [left] {\y};		
        \end{tikzpicture}}
      \caption{Graphical representation of the ground truth of our simulated
        experiments. \emph{Left:} The entries from the diagonal matrix $X$ with
        their corresponding indices and values. Only the indicated entries are
        non-zero. \emph{Right:} The corresponding spatial distribution.}
      \label{fig:plotsol}
\end{figure}
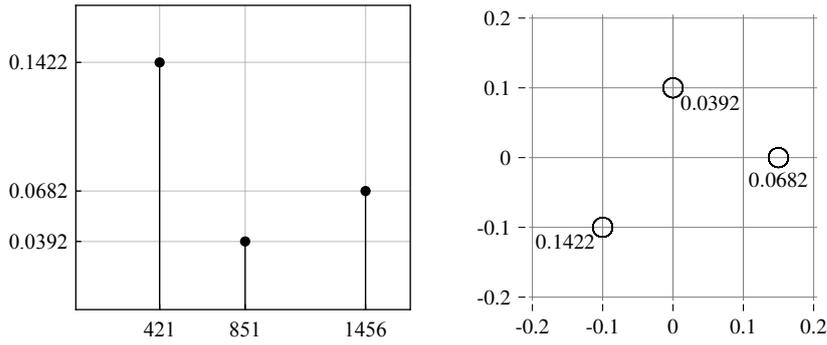

Our tests use different data setups. The considered data stems from simulated
experiments and therefore, an exact ground truth is to our avail. A detailed
presentation of the considered setup is given in Section~\ref{sec:basic-setup}.
In our setup we have $A\in\Complex^{64,1681}$ and $C\in\Complex^{64,64}$, such
that $X$ must be a diagonal matrix in $\Complex^{1681,1681}$. The exact solution
is supposed to contain 3 sources with different strengths. The corresponding
diagonal matrix $X_{\text{sol}}$ has a non-zero entry in the positions $421$,
$851$ and $1456$ of its main diagonal. The corresponding values are $0.1422$,
$0.0392$ and $0.0682$ respectively. A plot visualising the solution in form of
the diagonal entries from the matrix $X$ as well as their spatial distribution
is given in Fig.~\ref{fig:plotsol}.\par
Our numerical tests cover the commonly occurring use-cases where different
frequencies and noise corruption must be handled. Data corruption can occur as
noise applied onto the recorded signals or as inaccuracies in the measurement of
the position. The latter case is much more difficult to handle since there is no
clear way to deduce the correct location from faulty data. This problem is
severely ill-posed. We remark, that the physical setup is encoded in the matrix
$C$. The matrix $A$ is given by the sound propagation model detailed in
Section~\ref{sec:Introduction}.\par
We will refrain from using any error measure, such as the popular $\ell_{2}$- or
$\ell_{1}$-distance since they do not allow to reflect the sparsity requirements
in a convincing manner. Neither of these distances reflect the concrete physical
properties of the underlying problem and further, they cannot distinguish
between an additional source with a certain strength and a wrongly estimated
source with the same deviation. In both cases the distance to the ground truth
would be the same. Therefore, we will only provide graphical representations of
the solution. The plots accurately visualise the determined number of sources.
Furthermore, since our setups are very sparse, it is easy to compare the
strengths of the individual sources to the ground truth.\par
In the following we evaluate our model from \eqref{eq:11} and from \eqref{eq:12}
and compare our results to reference solutions obtained with well established
methods from the literature. We do not evaluate the model from \eqref{eq:10}
explicitly, since it is a special case of the formulation in \eqref{eq:11} where
all entries in the weighting matrix take the same value. Furthermore,
\eqref{eq:10} provides no mechanism to favour a diagonal structure of the
solution $X$. Therefore, we expect it to perform significantly worse than the
other two strategies in the experiments considered in this work. Nevertheless,
we point out that \eqref{eq:10} is probably a valuable and competitive strategy
when the sources are not assumed to be uncorrelated. In this case non-zero
entries outside of the main diagonal are to be expected.\par
%
% ------------------------------------------------------------------------------
\subsection{Basic Setup}
\label{sec:basic-setup}
% ------------------------------------------------------------------------------
% 
The exemplary cross-spectral matrix $C$ of microphone signals used for this
study is synthesised from simulated data, following a typical data processing
used with microphone array measurements for evaluation in the frequency
domain.\par
\begin{figure}
      \small
      \caption{\emph{Left:} Data generation parameters. \emph{Right:}
        Distribution of the 64 microphones (filled circles) with the simulated
        sources (crossed circles). Also shown is the discretised grid of focus
        points. The array is positioned $\SI{0.3}{\metre}$ above the sources,
        parallel to the focus grid.}
      \label{tab:datagen}
      \begin{minipage}{0.58\linewidth}
            \normalsize\vspace*{1em}
            \begin{tabular}{ll}
              \toprule \addlinespace
              \multirow{3}{2.5cm}{Source positions and $c_{\text{RMS},\SI{1}{\metre}}$}
                                                     & $(-0.1, -0.1, 0.3)$\;m, $\SI{1}{\pascal}$  \\
                                                     & $(0.15, 0.0, 0.3)$\,m, $\SI{0.7}{\pascal}$ \\
                                                     & $(0.0, 0.1, 0.3)$\,m, $\SI{0.5}{\pascal}$  \\
              \addlinespace
              Microphones                            & $64$                                       \\
              Meas.\ time                            & $\SI{80}{\second}$                         \\
              Sampling rate                          & $\SI{51.2}{\kilo\hertz}$                   \\
              \addlinespace
              \acrshort{FFT} block size              & $128$ samples                              \\
              \multirow{2}{*}{\acrshort{FFT} window} & von Hann,                                  \\
                                                     & $50\%$ overlap                             \\
              \multirow{2}{*}{Focus grid}            & $41\times 41$ points                       \\
                                                     & $\Delta x = \Delta y = \SI{0.01}{\metre}$  \\
              \addlinespace
              \bottomrule
            \end{tabular}
      \end{minipage}
      \hfill
      \begin{minipage}{0.41\linewidth}
            \includegraphics[width=\linewidth]{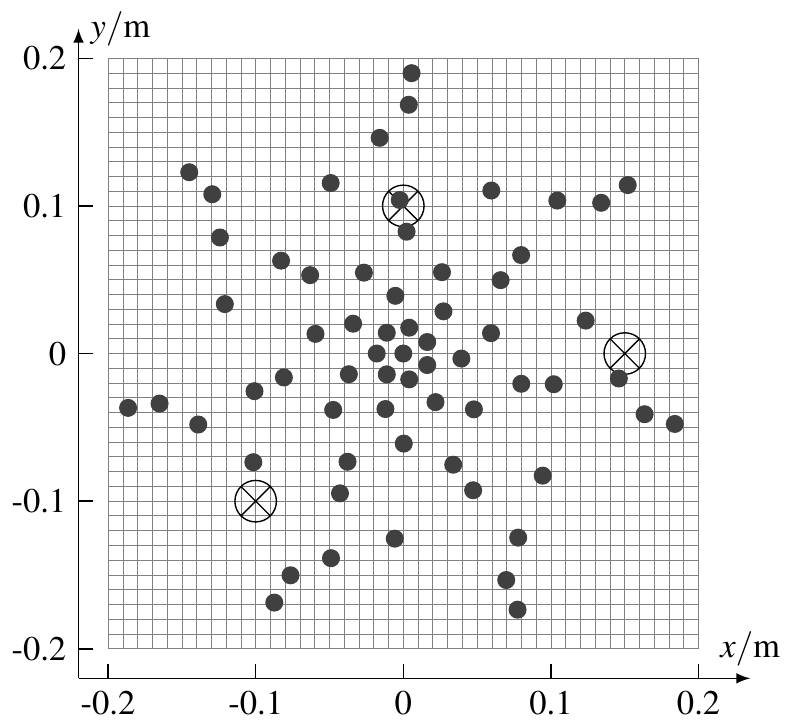} 
      \end{minipage}
\end{figure}

First, microphone time signals are simulated by calculating the sound pressures
caused by three monopole sources distributed at different positions in one plane
parallel to the array. These sources emit uncorrelated Gaussian white noise with
differing source strengths, defined by the \gls{RMS} value of the sound pressure
in $\SI{1}{\metre}$ distance from the respective source. Data generation
parameters and the array geometry with the relative source positions are
specified in Fig.~\ref{tab:datagen}.
Following Welch's method \cite{W1967}, the time signals are cut into $K$
overlapping blocks of equal length, onto which a \gls{FFT} is applied. For each
block, the cross-spectra between the microphone channels are calculated from the
vector of complex sound pressures $c_k$. Finally, the cross-spectral matrix is
estimated by averaging the cross-spectra:
\begin{equation}
\label{eq:CSM}
C = \frac{1}{K}\sum_{k}{c_k\,{c_k}^\top}~,\quad k=1, \ldots, K
\end{equation}
The assumed possible source locations form a square grid of regularly spaced
focus points, contain the actual source positions, and are arranged in a plane
parallel to the array. The extension of the focus grid is of the same order of
magnitude as that of the array (see Fig.~\ref{tab:datagen}). With a sampling
rate of $\SI{51.2}{\kilo\hertz}$ and a \gls{FFT} block size of 128 we can compute
63 cross-spectral matrices from the microphone signals. Each matrix represents a
discrete frequency band with band centre frequencies between $\SI{400}{\hertz}$
and $\SI{25.2}{\kilo\hertz}$ and with a band width of $\SI{400}{\hertz}$. All
experiments in this paper use a cross-spectral matrix with a band centre
frequency of $\SI{19.2}{\kilo\hertz}$.\par
The noisy data is generated by corrupting the microphone signals directly. Each
microphone is corrupted individually. The noise applied to the signals is
uncorrelated and has a relative strength of $3.33$, $4.76$ and $6.25$
respectively, when compared to the signal strength. In particular, the noise is
stronger than the actual signal.\par
Additionally, a data set with erroneous microphone positions was generated. For
this, the nominal position of each microphone was disturbed randomly within the
array plane. The average deviation from the true position is
$\SI{0.009}{\meter}$. The array has an aperture of $\SI{0.4}{\meter}$.
% 
% ------------------------------------------------------------------------------
\subsection{Using the Weighted Model from Equation~\eqref{eq:11}}
\label{sec:NoStructWithWeight}
% ------------------------------------------------------------------------------
%
In this section we analyse the performance of our model stated in \eqref{eq:11}.
Our first test uses optimal uncorrupted data. We set the number of Bregman
iterations, alternating minimisations and gradient descent steps to
$75$/$3$/$20$. The weighting matrix has the value $1$ on its main diagonal and
$10^{6}$ in all other positions. The Bregman regularisation weight has the value
$10^{5}$. With these parameters the run time of a pure Matlab implementation was
about $20$ minutes for a single frequency band. Our algorithms are always
initialised with a solution having the value $1$ at all positions. The run time
of the post processing is negligible. The k-means algorithm always converged
within a few iterations and had run times below one second.
%
% ------------------------------------------------------------------------------
\subsubsection{Using Perfect Data}
\label{sec:using-noise-free}
\begin{figure}
  \centering
  \includegraphics{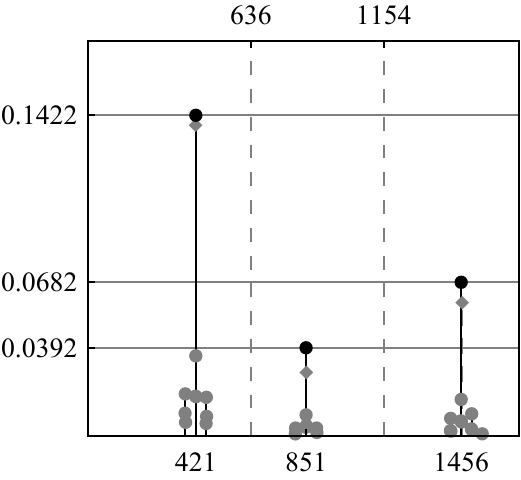}
  \hspace*{2em}
  \raisebox{0.05\height}{\includegraphics{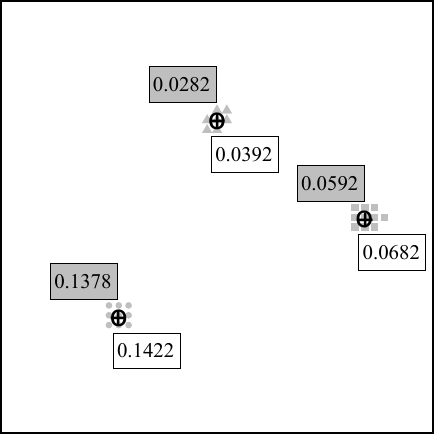}}
  \caption{Results obtained from Algorithm~\ref{alg:SplitBregmanNoStructure} for
    solving \eqref{eq:11} for the noise-free setup. \emph{Left:} Diagonal entries
    returned by the algorithm. The grey discs represent the solution without
    post processing, whereas the grey squares denote the post-processed
    solution. The dashed vertical lines indicate the cluster boundaries obtained
    by our k-means algorithm. We see three distinct clusters of entries centred
    around the exact solution (marked in black). \emph{Right:} The corresponding
    graphical representation on the 2D Grid. The filled grey shapes indicate the
    positions of the non-zero entries from the main diagonal. Each shape
    represents one of the clusters returned by our post-processing strategy. The
    crosses and the grey labels mark the position and strength of the computed
    centres of the clusters. The white labels and circles depict the ground
    truth solution.}
  \label{fig:f48e1}
\end{figure}

% ------------------------------------------------------------------------------
%
The results obtained from the perfect data are visualised in
Fig.~\ref{fig:f48e1}. As we can see, the reconstruction is almost perfect. The
diagonal has very few non-zero entries, that cluster around the exact solutions.
All the other entries are smaller than $10^{-5}$ in magnitude and have been
excluded from the plot to ease the visualisation. After the post-processing step
with our k-means approach we obtain the reconstruction in the right plot in
Fig.~\ref{fig:f48e1}.
%
% ------------------------------------------------------------------------------
\subsubsection{Using Noisy Data}
\label{sec:using-noise}
\begin{figure}
      \centering
      \includegraphics{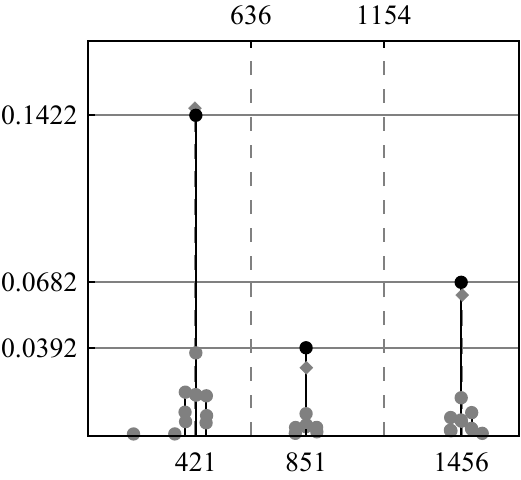}
      \hspace*{2em}
      \raisebox{0.05\height}{\includegraphics{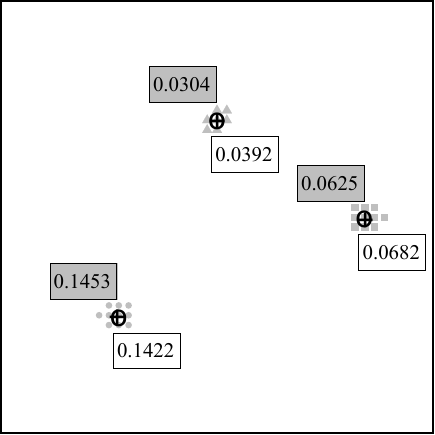}}
      \caption{Results obtained from Algorithm~\ref{alg:SplitBregmanNoStructure}
        for solving \eqref{eq:11} for the noisy setup. \emph{Left:} Diagonal
        entries returned by the algorithm. The grey discs represent the solution
        without post processing, whereas the grey squares denote the
        post-processed solution. The dashed vertical lines indicate the cluster
        boundaries obtained by our k-means algorithm. We see three distinct
        clusters of entries centred around the exact solution (marked in black).
        \emph{Right:} The corresponding graphical representation on the 2D Grid.
        The filled grey shapes indicate the positions of the non-zero entries
        from the main diagonal. Each shape represents one of the clusters
        returned by our post-processing strategy. The crosses and the grey
        labels mark the position and strength of the computed centres of the
        clusters. The white labels and circles depict the ground truth solution.
        As we can see, our algorithm behaves well, even in the presence of
        noise. The results are almost identical to those from
        Fig.~\ref{fig:f48e1}}
      \label{fig:f48n1}
\end{figure}

% ------------------------------------------------------------------------------
%
Similar findings to those from the previous paragraph are also obtainable when
the data is corrupted by noise. The results are given in Fig.~\ref{fig:f48n1}.
They show that our approach is robust to noise. The obtained solutions,
especially after the post-processing, hardly differ from those from the perfect
setup.
%
% ------------------------------------------------------------------------------
\subsubsection{Using Data with Positional Errors}
\label{sec:using-pos-error}
\begin{figure}
      \centering
      \includegraphics{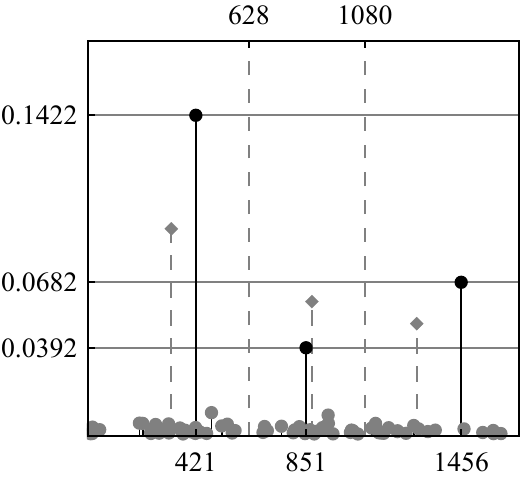}
      \hspace*{2em}
      \raisebox{0.05\height}{\includegraphics{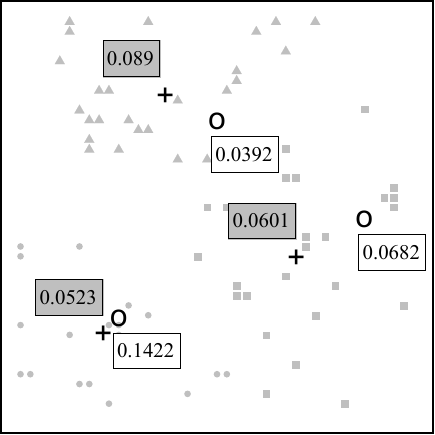}}
      \caption{Graphical representation of the diagonal matrix using data with
        erroneous positional information as returned by
        Algorithm~\ref{alg:SplitBregmanNoStructure}. \emph{Left:} A
        straightforward application yields an undesirable outcome. There are
        many non-zero entries along the main diagonal. The grey discs represent
        the solution without post processing, whereas the grey squares denote
        the post-processed solution. The dashed vertical lines indicate the
        cluster boundaries obtained by our k-means algorithm. After applying the
        post processing steps we obtain a convincing localisation. The black
        circles denote the exact ground truth. \emph{Right:} The corresponding
        graphical representation on the 2D Grid. The filled grey shapes indicate
        the positions of the non-zero entries from the main diagonal. Each shape
        represents one of the clusters returned by our post-processing strategy.
        The crosses and the grey labels mark the position and strength of the
        computed centres of the clusters. The white labels and circles depict
        the ground truth solution.}
      \label{fig:f48p1}
\end{figure}

% ------------------------------------------------------------------------------
%
Finally, the data with positional errors yields the results presented in
Fig.~\ref{fig:f48p1}. Without the k-means post processing the solution is hardly
usable. However, with the additional clustering step we obtain a rather good
reconstruction. Let us emphasise that the handling of the data sets with errors
in the positional data is a severely ill posed problem for which no satisfactory
solution exists in the literature so far.
%
% ------------------------------------------------------------------------------
\subsection{Using the Model with Structural Constraints from
  Equation~\eqref{eq:12}}
\label{sec:StructNoWeight}
% ------------------------------------------------------------------------------
%
In this section we want to demonstrate the benefits of integrating structural
constraints into the model and optimisation process. We test our
Algorithm~\ref{alg:SplitBregmanStructure} on the same data as in the previous
section. Applying Algorithm~\ref{alg:SplitBregmanStructure} with $75$ outer
iterations, $4$ alternating minimisation steps and $10$ gradient descent steps.
%
% ------------------------------------------------------------------------------
\subsubsection{Using Perfect Data}
\label{sec:using-noise-free-2}
\begin{figure}
      \centering
      \includegraphics{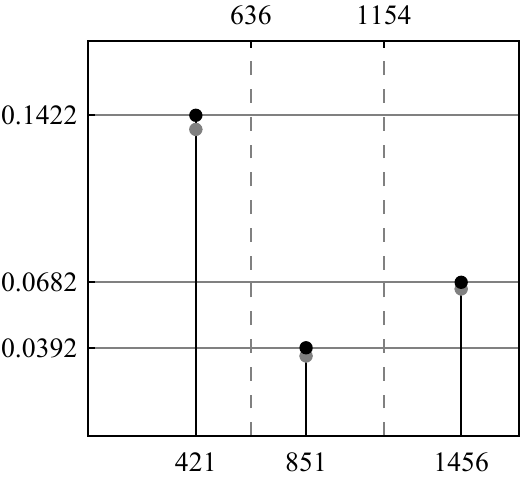}
      \hspace*{2em}
      \raisebox{0.05\height}{\includegraphics{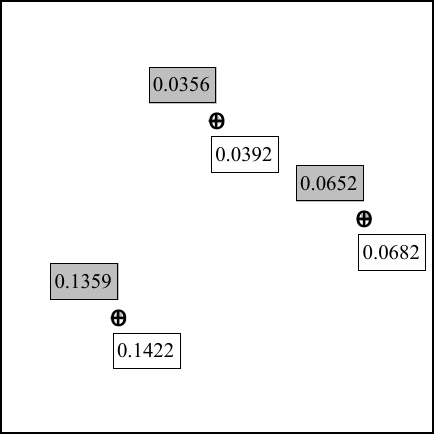}}
      \caption{Results obtained from Algorithm~\ref{alg:SplitBregmanStructure}
        for solving \eqref{eq:12} for the noise-free setup. \emph{Left:} Diagonal
        entries returned by the algorithm. The grey discs represent the solution
        without post processing, whereas the grey squares denote the
        post-processed solution. The dashed vertical lines indicate the cluster
        boundaries obtained by our k-means algorithm. Even without the post
        processing the solution is nearly indistinguishable from the exact
        solution (marked in black). \emph{Right:} The corresponding graphical
        representation on the 2D Grid. The filled grey shapes indicate the
        positions of the non-zero entries from the main diagonal. Each shape
        represents one of the clusters returned by our post-processing strategy.
        The crosses and the grey labels mark the position and strength of the
        computed centres of the clusters. The white labels and circles depict
        the ground truth solution.}
      \label{fig:SolT20re}
\end{figure}

% ------------------------------------------------------------------------------
%
The results obtained from the perfect data are visualised in
Fig.~\ref{fig:SolT20re}. Here, the sparsity weight was set to $10$ and the
Bregman regularisation parameter was set to $10^{4}$. As we can see, the
obtained result is almost identical to the reference ground truth.\par
%
% ------------------------------------------------------------------------------
\subsubsection{Using Noisy Data}
\label{sec:using-noise-2}
\begin{figure}
      \centering
      \includegraphics{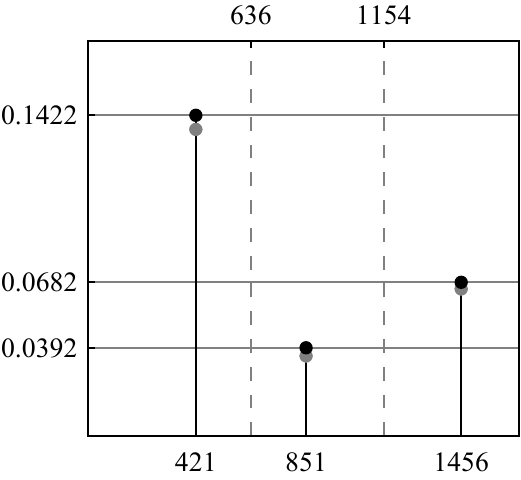}
      \hspace*{2em}
      \raisebox{0.05\height}{\includegraphics{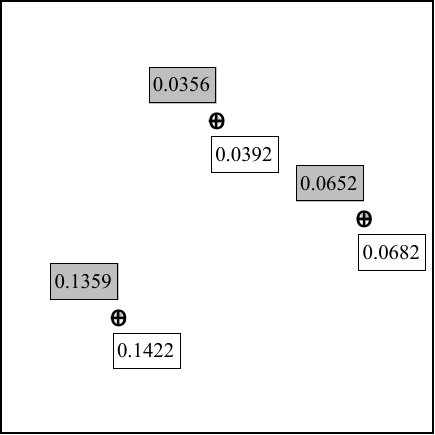}}
      \caption{Results obtained from Algorithm~\ref{alg:SplitBregmanStructure}
        for solving \eqref{eq:12} for the noisy data setup. Again, our approach
        shows a very favourable behaviour in presence of noise. Our results are
        identical to those from the noise free test case.}
      \label{fig:d48n}
\end{figure}

% ------------------------------------------------------------------------------
%
Figure~\ref{fig:d48n} depicts our reconstruction from the noise corrupted
signal. As we can see, our approach is robust and still yields almost ideal
results. The parameters were exactly the same as in the noise free case. Our
experiments suggest that the method works reliably, for almost any reasonable
parameter setting. We consider this a further asset of our strategy.\par
%
% ------------------------------------------------------------------------------
\subsubsection{Using Data with Positional Errors}
\label{sec:using-pos-error-2}
\begin{figure}
      \centering
      \includegraphics{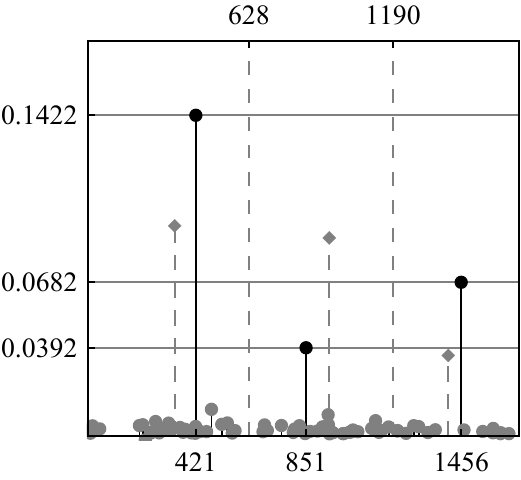}
      \hspace*{2em}
      \raisebox{0.05\height}{\includegraphics{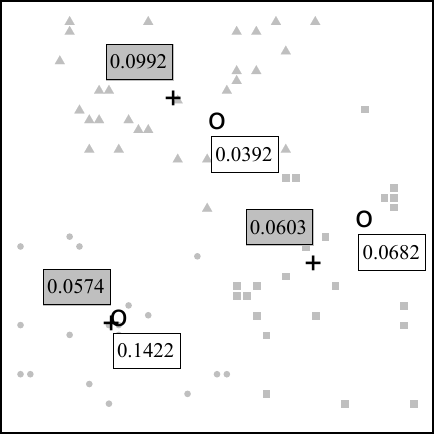}}
      \caption{Graphical representation of the of the diagonal matrix using data
        with erroneous positional information as returned by
        Algorithm~\ref{alg:SplitBregmanStructure}. \emph{Left:} A
        straightforward application yields a undesirable outcome. There are many
        non-zero entries along the main diagonal. Also here, the post processing
        improves the results significantly. \emph{Right:} The corresponding
        graphical representation on the 2D Grid. The filled grey shapes indicate
        the positions of the non-zero entries from the main diagonal. Each shape
        represents one of the clusters returned by our post-processing strategy.
        The crosses and the grey labels mark the position and strength of the
        computed centres of the clusters. The white labels and circles depict
        the ground truth solution.}
      \label{fig:d48p1}
\end{figure}

% ------------------------------------------------------------------------------
%
The data set with the encoded erroneous positional information has proven to be
the most difficult to handle. Executing our method with $75$ Bregman iterations,
$6$ alternating minimisations and $50$ gradient descent steps yields the results
shown in Fig.~\ref{fig:d48p1}. Other parameter choices did not yield
significantly better results. As we can see, the entries on the main diagonal of
the matrix are not as sparse as before and they are off by almost a factor $10$.
In order to remedy the situation, we apply our k-means based post-processing.
The resulting clusters, and their centroid positions are presented in
Fig.~\ref{fig:d48p1}.
The imaginary part of every entry is smaller than $10^{-5}$ in magnitude. Also,
all entries smaller in absolute value than $10^{-3}$ have been omitted to
improve the visualisation. Clearly, our result consists of 3 distinct clusters
which are well localised. The position of each cluster coincides with the
coordinates of one peak from the ground truth. However, the scaling does not
match. Yet, it is interesting to note that the sum of all entries in a single
cluster almost gives the sought solution value.
%
% ------------------------------------------------------------------------------
\subsection{Comparison to Other Methods from the Literature}
\label{sec:Comparison}
% ------------------------------------------------------------------------------
% 
Let us shortly compare the output of our schemes to well established methods
from the literature. To this end we have evaluated our experimental data with
the \gls{CLEANSC} \cite{S2007} and \gls{CMF} \cite{BE2004, YLSC2008} strategies.
\gls{CLEANSC} is an extension of the venerable CLEAN algorithm for deconvolution
\cite{H1974}. The \gls{CMF} method is similar in spirit as our optimisation
models from \eqref{eq:10}. We use a soft constraint on the sparsity, whereas
\gls{CMF} uses a hard constraint.\par
The methods have been implemented and used ``as-is'' without any additional post
or preprocessing of the data. We remark that the considered methods are tailored
towards yielding sparse solutions. As a consequence, post processing steps, such
as our clustering approach, make little sense in this context.
%
% ------------------------------------------------------------------------------
\subsubsection{Evaluation of the Covariance Matrix Fitting Method}
\label{sec:cmf}
\begin{figure}
      \centering
      \includegraphics{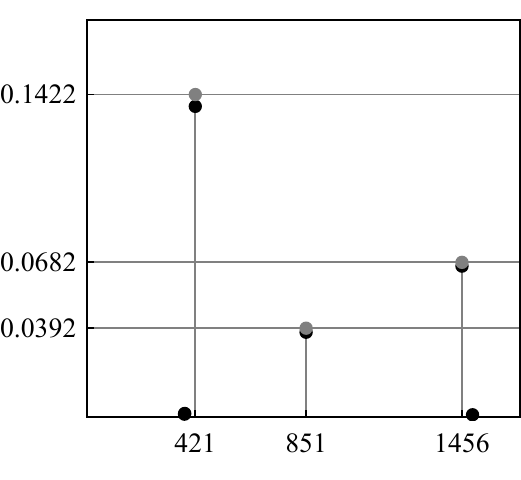}
      \hspace*{2em}
      \raisebox{0.05\height}{\includegraphics{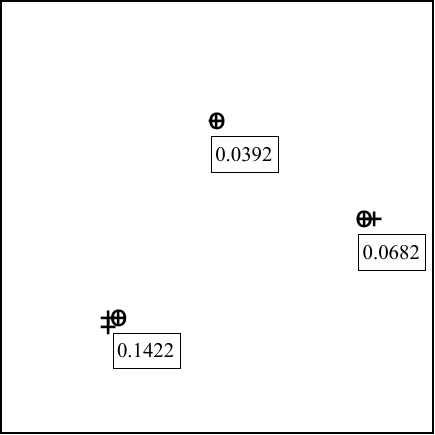}}
      \caption{Solution from perfect data with the \gls{CMF} method. Except for
        two very weak additional sources detected in the vicinity of the
        strongest sources the method finds a near perfect solution. The
        positions of the three strongest sources coincide with the ground truth,
        only the intensities are slightly off. \emph{Left:} Graphical
        representation of the main diagonal. \emph{Right:} Corresponding spatial
        distribution. Crosses mark the obtained solution whereas the labelled
        circles denote the ground truth. The labels of the obtained solution
        have been omitted to avoid unnecessary cluttering of the visualisation.}
      \label{fig:cmfperfect}
\end{figure}

\begin{figure}
      \centering
      \includegraphics{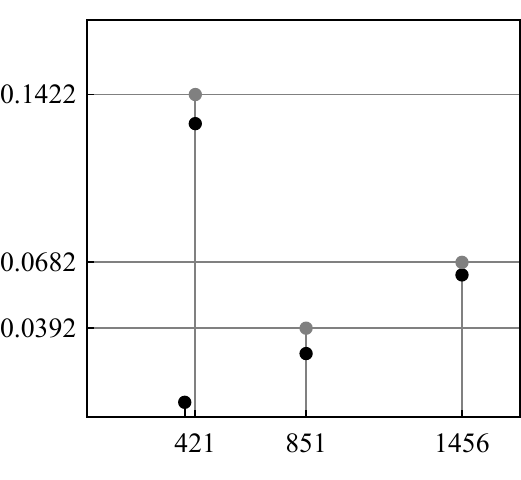}
      \hspace*{2em}
      \raisebox{0.05\height}{\includegraphics{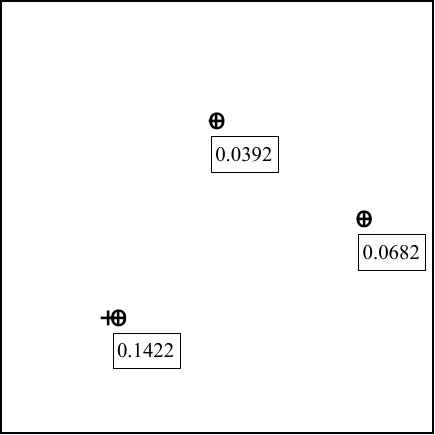}}
      \caption{Solution from noisy data with the \gls{CMF} method. Except for
        two very weak additional sources detected in the vicinity of the
        strongest sources the method finds a near perfect solution. The solution
        of the \gls{CMF} algorithm is almost unaffected by the addition of
        noise. \emph{Left:} Graphical representation of the main diagonal.
        \emph{Right:} Corresponding spatial distribution. Crosses mark the
        obtained solution whereas the labelled circles denote the ground truth.
        The labels of the obtained solution have been omitted to avoid
        unnecessary cluttering of the visualisation.}
      \label{fig:cmfnoisy}
\end{figure}

\begin{figure}
      \centering
      \includegraphics{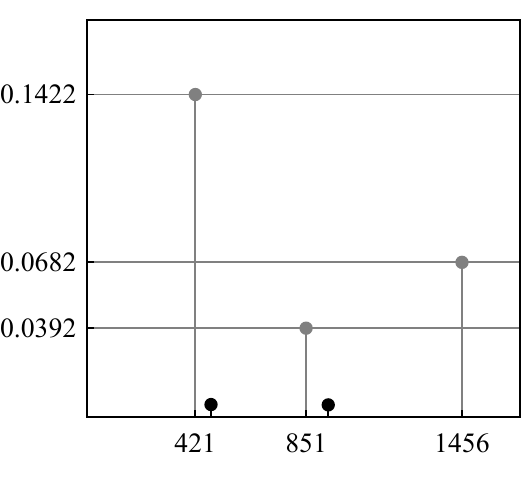}
      \hspace*{2em}
      \raisebox{0.05\height}{\includegraphics{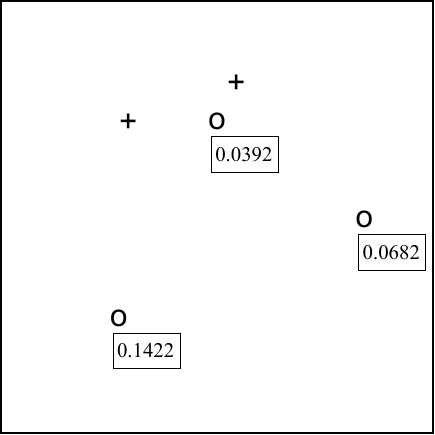}}
      \caption{Solution from data with erroneous positions with the \gls{CMF}
        Method. In this case the algorithm breaks completely down. Not only the
        positions and number of sources are wrong, but also the signal strengths
        are in a completely wrong range. \emph{Left:} Graphical representation
        of the main diagonal. \emph{Right:} Corresponding spatial distribution.
        Crosses mark the obtained solution whereas the labelled circles denote
        the ground truth. The labels of the obtained solution have been omitted
        to avoid unnecessary cluttering of the visualisation.}
      \label{fig:cmfpos}
\end{figure}

% ------------------------------------------------------------------------------
%
Using the \gls{CMF} method we obtain the results presented in
Figs.~\ref{fig:cmfperfect}, \ref{fig:cmfnoisy} and \ref{fig:cmfpos}. The figures
represent the use-cases with perfect and noisy data as well as with erroneously
encoded positions respectively.
As we can observe in Fig.~\ref{fig:cmfperfect}, the \gls{CMF} method fails to
accurately detect the correct number of sources. It adds two additional sources
(with almost negligible strength). Nevertheless, the three largest peaks
coincide almost perfectly with the ground truth solution from
Fig.~\ref{fig:plotsol}.
In the case of noisy data in Fig.~\ref{fig:cmfnoisy}, the \gls{CMF} approach
works equally well. Position and magnitude of the three largest peaks are very
accurate.
Finally, the method completely fails if the positions are erroneously encoded
(c.f.\ Fig.~\ref{fig:cmfpos}). Not only the position, but also the number and
the strength of the signals are incorrectly detected.
% 
% ------------------------------------------------------------------------------
\subsubsection{Evaluation of the Clean-SC Method}
\label{sec:csc}
\begin{figure}
      \centering
      \includegraphics{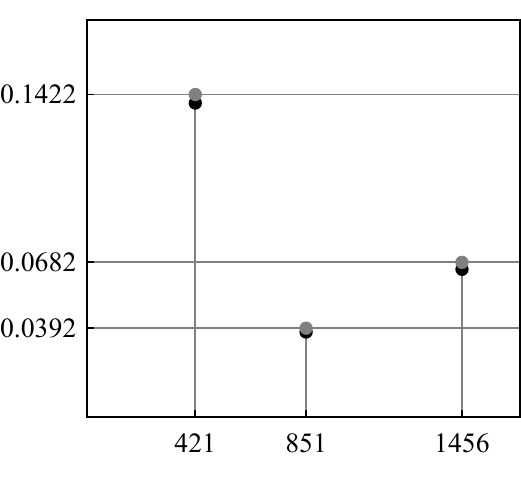}
      \hspace*{2em}
      \raisebox{0.05\height}{\includegraphics{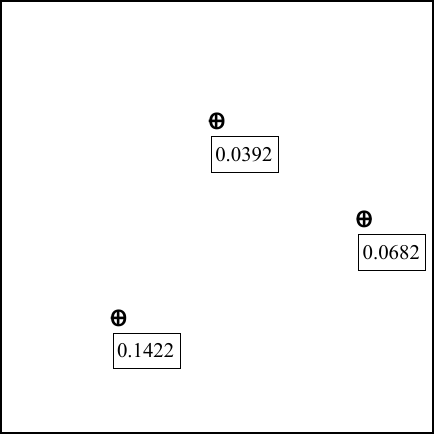}}
      \caption{Solution from perfect data with the \gls{CLEANSC} Method. The
        number of source, their position and their strength are detected with a
        very high accuracy. \emph{Left:} Graphical representation of the main
        diagonal. \emph{Right:} Corresponding spatial distribution. Crosses mark the
        obtained solution whereas the labelled circles denote the ground truth.
        The labels of the obtained solution have been omitted to avoid
        unnecessary cluttering of the visualisation.}
      \label{fig:cscperfect}
\end{figure}

\begin{figure}
      \centering
      \includegraphics{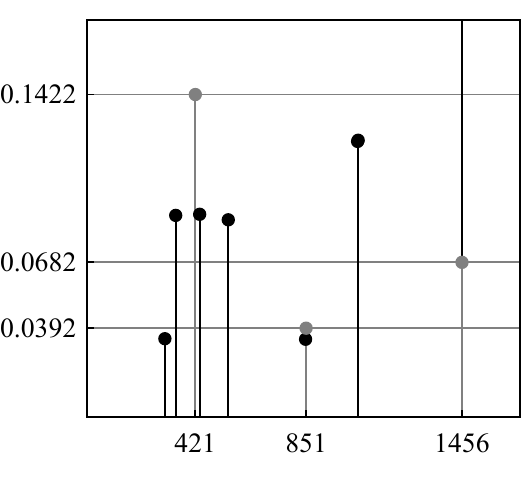}
      \hspace*{2em}
      \raisebox{0.05\height}{\includegraphics{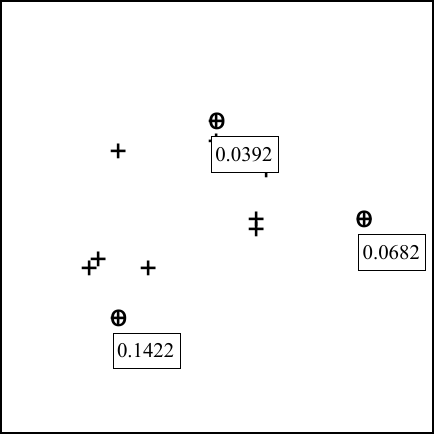}}
      \caption{Solution from noisy data with the \gls{CLEANSC} Method. The
        algorithm suffers more from noisy data than \gls{CMF}. The number of
        sources is wrongly estimated and the magnitude of the two largest peaks
        is completely wrong. The ticks on the y axis denote the correct
        magnitude. As we can see, the third peak is about 3 times larger than it
        should. \emph{Left:} Graphical representation of the main diagonal.
        \emph{Right:} Corresponding spatial distribution. Crosses mark the
        obtained solution whereas the labelled circles denote the ground truth.
        The labels of the obtained solution have been omitted to avoid
        unnecessary cluttering of the visualisation.}
      \label{fig:cscnoisy}
\end{figure}

\begin{figure}
      \centering
      \includegraphics{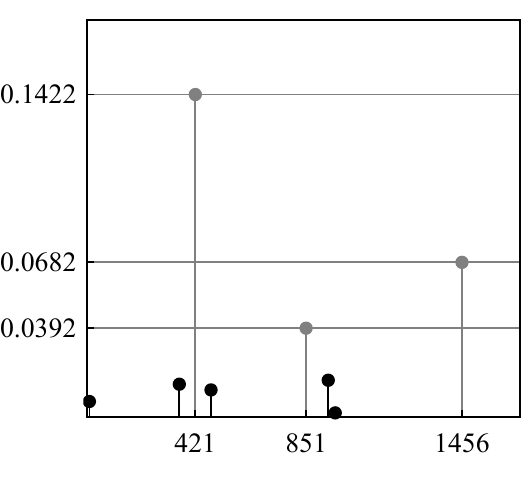}
      \hspace*{2em}
      \raisebox{0.05\height}{\includegraphics{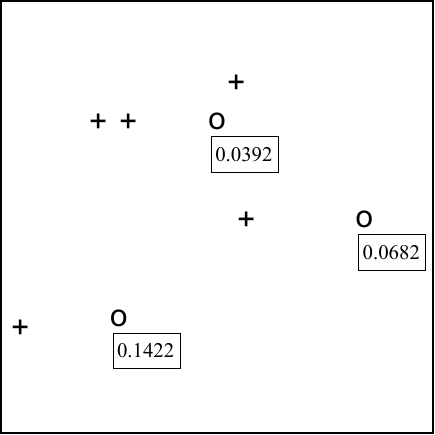}}
      \caption{Solution from data with erroneous positions with the
        \gls{CLEANSC} Method. \gls{CLEANSC} behaves similarly as \gls{CMF}.
        \emph{Left:} Graphical representation of the main diagonal.
        \emph{Right:} Corresponding spatial distribution. Crosses mark the
        obtained solution whereas the labelled circles denote the ground truth.
        The labels of the obtained solution have been omitted to avoid
        unnecessary cluttering of the visualisation.}
      \label{fig:cscpos}
\end{figure}

% ------------------------------------------------------------------------------
%
Using the \gls{CLEANSC} method we obtain the results presented in
Figs.~\ref{fig:cscperfect}, \ref{fig:cscnoisy} and \ref{fig:cscpos}. The figures
represent the use-cases with perfect and noisy data, as well as with erroneously
encoded positions respectively.
In presence of perfect data, the \gls{CLEANSC} method works extraordinarily well
and even outperforms the \gls{CMF} algorithm. The yielded results are nearly
perfect, see Fig.~\ref{fig:cscperfect}.
On the other hand, if the data is corrupted by noise, then the \gls{CLEANSC} strategy
cannot compete with \gls{CMF}. The method yields too many sources as can be seen in
Fig.~\ref{fig:cscnoisy}. Furthermore, the strength of certain sources is
completely wrong. The largest peak in Fig.~\ref{fig:cscnoisy} is nearly three
times larger than it should.
In presence of positional errors the method performs similarly to \gls{CMF},
the findings are little convincing, as can be observed in Fig.~\ref{fig:cscpos}.
%
% ------------------------------------------------------------------------------
\subsection{Conclusion}
\label{sec:EvalConclusion}
% ------------------------------------------------------------------------------
%
The following conclusions can be drawn from the evaluation.
\begin{enumerate}
\item State-of-the-art methods can handle perfect data as well as noisy data
      quite well.
\item Our algorithms yield competitive results for perfect and noisy data. One
      advantage of our method is that the number of sources is always accurate.
\item \gls{CMF} and \gls{CLEANSC} cannot handle data sets with wrongly encoded
      positions at all.
\item Our algorithms do not yield perfect solutions for wrongly encoded
      positions either, nevertheless the quality is by far superior to the
      findings from \gls{CMF} and \gls{CLEANSC}.
\end{enumerate}
In the total, our new algorithms provide accurate and more robust alternatives
compared to actual state-of-the-art methods.
%
% ------------------------------------------------------------------------------
\section{Summary and Outlook}
\label{sec:ConclusionOutlook}
% ------------------------------------------------------------------------------
%
We have presented three approaches to retrieve a sparse set of sound source
locations from acoustic measurements. Our novel algorithms apply findings from
matrix differential calculus to sparsity favouring convex optimisation models,
which we have applied onto complex matrix valued settings. Our efforts to
integrate the structural constraints, such as the shape of the matrix to be
optimised into the optimisation problem, have been rewarded by very convincing
results. Not only did we benefit from the reduced memory footprint and faster
run times, but also the obtained results are much more accurate than in the
unconstrained framework. Empirical studies show that our algorithms are
competitive to state-of-the-art methods from the literature for perfect and
noisy data. In presence of ill-posed setups, we benefit from our post processing
strategies to obtain reasonable results.\par
In the future, we would like to provide further improvements to the accuracy of
the results. Especially the handling of inaccuracies in the recording of the
positions will be of interest to us. Further investigations on acceleration
methods for Bregman schemes will also be a topic of our research. We are
confident that a significant reduction of the run time is possible.
%
% \begin{acknowledgements}
%   If you'd like to thank anyone, place your comments here
%   and remove the percent signs.
% \end{acknowledgements} 
% 
% ------------------------------------------------------------------------------
\bibliographystyle{spmpsci}
\bibliography{references}
% ------------------------------------------------------------------------------
\end{document}